\documentclass[]{article}
\usepackage{graphicx}
\usepackage{amsfonts}
\usepackage{amsmath}
\usepackage{amssymb}
\usepackage{fancyhdr}
\usepackage{titlesec}
\usepackage{indentfirst}
\usepackage{booktabs}
\usepackage{verbatim}
\usepackage{color}
\usepackage{amsthm}
\usepackage{caption}
\usepackage{subcaption}
\usepackage{physics}
\usepackage{mathtools}
\usepackage{tikz}
\usepackage{enumerate}
\usepackage{stmaryrd}

\usepackage{hyperref}
\usepackage{cleveref}
\usepackage[page,header]{appendix}
\usepackage{titletoc}

\newcommand{\rd}{\,\mathrm{d}}

\numberwithin{equation}{section}
\newtheorem{theorem}{Theorem}[section]

\newtheorem{remark}[theorem]{Remark}

\newcommand{\rom}[1]{\text{\MakeUppercase{\romannumeral #1}}}
\newcommand{\rombracket}[1]{(\text{\MakeUppercase{\romannumeral #1}})}

\renewcommand{\bf}{\mathbf{f}}

\newcommand{\bF}{\mathbf{F}}
\newcommand{\bv}{\mathbf{v}}
\newcommand{\e}{\mathrm{e}}
\newcommand{\bu}{\mathbf{u}}
\newcommand{\bI}{\mathbf{I}}
\newcommand{\bx}{\mathbf{x}}
\newcommand{\bU}{\mathbf{U}}
\newcommand{\dt}{\Delta t}
\newcommand{\gradv}{\nabla_\bv}

\newcommand{\pdt}{\partial_t}

\newcommand{\bD}{\mathbf{D}}
\newcommand{\bQ}{\mathbf{Q}}
\newcommand{\bK}{\mathbf{K}}
\newcommand{\bM}{\mathbf{M}}
\newcommand{\up}{{\mathrm{up}}}
\renewcommand{\dv}{\Delta v}
\newcommand{\dx}{\Delta {x}}

\begin{document}

\title{Implicit Dynamical Tensor Train Approximation \\for Kinetic Equations with Stiff Fokker--Planck Collisions}


\author{
Geshuo Wang\footnote{Department of Applied Mathematics, University of Washington, Seattle, WA 98195 (geshuo@uw.edu).} \
\ and \ 
Jingwei Hu\footnote{Department of Applied Mathematics, University of Washington, Seattle, WA 98195 (hujw@uw.edu). Corresponding author.}
}

\maketitle

\begin{abstract}
Low-rank methods for kinetic equations have attracted increasing attention due to their effectiveness in reducing the high dimensionality of phase space. In our previous work [G.~Wang \& J.~Hu, \emph{J.~Comput.~Phys.} 558 (2026) 114884], we developed a dynamical low-rank method based on the projector-splitting integrator in tensor-train (TT) format, in which explicit time integration is employed in all substeps. 
As a result, the method is subject to severe stability constraints in the strongly collisional regimes. In this paper, we consider kinetic equations with the (nonlinear) Fokker--Planck collision operator and develop a dynamical low-rank method that employs implicit or implicit-explicit (IMEX) discretizations in appropriate substeps to overcome stiffness. In these implicit substeps, the resulting equations can be formulated as matrix or tensor Sylvester equations, for which we propose efficient direct solvers by exploiting their underlying structure. 
The overall computational cost of the proposed method scales linearly with respect to the number of grid points in a single velocity dimension, comparable to that of a fully explicit low-rank scheme.
We demonstrate the accuracy and efficiency of the proposed method on several representative kinetic test problems.
\end{abstract}

\section{Introduction}
Kinetic equations describe the evolution of probability density distributions in phase space and 
play a fundamental role in many areas such as rarefied gas dynamics \cite{Cercignani} and plasma physics \cite{Villani02}.
A major challenge in their numerical simulation is the high dimensionality of phase space, which leads to prohibitive memory and computational costs for conventional discretization methods. 
In addition, many physically relevant regimes involve collision operators, such as the Fokker--Planck operator, which may introduce stiffness and further restrict the choice of numerical schemes.

To address the curse of dimensionality, low-rank methods have emerged as an effective approach by exploiting the intrinsic structure of the distribution function \cite{einkemmer2025review}, leading to significant reductions in both memory usage and computational cost. 
Various low-rank representations have been proposed for kinetic equations. 
A classical approach treats the distribution function as a function of space and velocity (e.g., \cite{einkemmer2018low,hu2022adaptive}) and applies matrix-based low-rank approximations via singular value decomposition (SVD). 
Alternatively, one may separate different directions in velocity space while treating the spatial variable as a parameter. 
This strategy is motivated by the observation that the equilibrium state, namely the Maxwellian distribution, is separable in the velocity variables. 
As a result, each spatial grid point can be associated with either a matrix-SVD representation \cite{galindoolarte2025nodal} or a tensor-based representation \cite{wang2026dynamical}.
For problems involving more than two variables, tensor-based representations are naturally employed. 
Among these, the Tucker or hierarchical Tucker decomposition \cite{Grasedyck10} and the tensor-train representation \cite{oseledets2011tensor} 
are particularly popular because of their favorable computational properties. Besides these representations, various strategies have been developed for solving time-dependent problems. 
The step-and-truncate (SAT) approach \cite{kormann2015semi,ehrlacher2017dynamical} advances the solution in time and subsequently applies truncation to control the rank. 
In contrast, the dynamical low-rank (DLR) approach \cite{koch2007dynamical} evolves the solution directly on the low-rank manifold by projecting the governing equation onto its tangent space. 
Although direct implementations of DLR may suffer from numerical instabilities associated with matrix inversions \cite{nonnenmacher2008dynamical,kieri2016discretized}, these difficulties can be effectively mitigated by projector-splitting integrators \cite{lubich2014projector,lubich2015time,lubich2018time}.
Another class of low-rank integrators that is also robust to small singular values is the augmented basis update \& Galerkin (BUG) integrator \cite{ceruti2022rank,ceruti2024parallel}.

Despite the advantages in terms of memory and computational cost, explicit low-rank schemes remain restricted by the Courant--Friedrichs--Lewy (CFL) condition. 
In regimes with moderate to strong collisions, for instance, the time step must be chosen to be very small to satisfy the stability constraints. The situation is even worse for Fokker--Planck-type collision operators, whose diffusive nature imposes a  parabolic CFL condition on explicit schemes. 
To overcome this limitation, implicit or implicit--explicit (IMEX) schemes have been incorporated into the low-rank framework. Earlier efforts include \cite{dolgov2012fast}, which uses alternating least squares (ALS) \cite{holtz2012alternating} and the density matrix renormalization group (DMRG) \cite{white1992density,white1993density,oseledets2011dmrg,oseledets2012solution} to solve the resulting nonlinear equations in low-rank format. Both approaches are iterative and require sweeping back and forth across the low-rank factors until convergence. Recently, a class of implicit low-rank methods has been proposed \cite{rodgers2023implicit,el2024krylov,nakao2025reduced,el2025Sylvester,appelo2025robust}. Although the details differ, all of these methods rely on iterative solvers in one way or another at every time step. In addition, truncation or rounding is generally required, since the rank of the solution typically increases during the implicit update process.


In contrast to previous works, we propose an efficient implicit dynamical low-rank method for solving the kinetic equations with stiff Fokker--Planck collisions. Our method is based on the projector-splitting integrator in tensor-train (TT) format \cite{wang2026dynamical} and has two distinct features: 
\begin{itemize}
\item The solution is evolved on a preselected fixed low-rank manifold in velocity space, so no truncation or rounding is needed. This representation becomes particularly efficient when the solution is close to the Maxwellian equilibrium. Rank adaptivity can be incorporated as needed, but it is not required for the implicit implementation.
\item The method does not require any iteration, meaning that the implicit step can be solved at essentially the same computational cost as an explicit low-rank method. This is achieved by exploiting the special structure of the Fokker--Planck operator and designing tailored direct solvers for the matrix and tensor Sylvester equations arising from the low-rank projection.
\end{itemize}


The rest of this paper is organized as follows. 
In \Cref{sec_integrator}, we first introduce the kinetic model and then discuss the discretization of the Fokker--Planck operator and its representation in TT format. The main part of this section presents the time integration of dynamical tensor trains using implicit or IMEX schemes, where we reduce the task to solving three special matrix or tensor Sylvester equations. In \Cref{sec_Sylvester_Solvers}, we develop efficient Sylvester solvers by exploiting the structure of the problem. The memory and computational costs are analyzed in \Cref{sec_cost}. 
Numerical results are presented in \Cref{sec_numerical_results}, followed by concluding remarks in \Cref{sec_conclusion}.

\section{Implicit dynamical tensor train approximation}
\label{sec_integrator}
Consider the kinetic equation of the form:
\begin{equation}
\label{eq_kinetic_equation}
    \pdt f(t,\bx,\bv) + \bv \cdot \nabla_\bx f(t,\bx,\bv) + \bF(t,\bx)\cdot \nabla_\bv f(t,\bx,\bv) = \eta \mathcal{Q}[f],
    \quad t>0 ,\quad \bx \in \Omega_\bx \subset \mathbb{R}^{3},
    \quad \bv \in \mathbb{R}^3,
\end{equation}
where $f$ is the distribution function depending on time $t$, position $\bx$, and velocity $\bv$;
$\bF$ is the acceleration due to external or self-consistent forces; and $\mathcal{Q}$ is the collision operator modeling interactions of particles with each other or with the surrounding environment. The parameter $\eta$ denotes the collision strength. 

In this paper, we assume $\mathcal{Q}$ to be the Fokker--Planck operator (often referred to as the Lenard--Bernstein \cite{LB58} or Dougherty \cite{Dougherty64} model in the plasma physics literature), which is widely used to describe Coulomb collisions between charged particles. It can be written in two equivalent forms as follows:
\begin{equation}
\label{eq_Fokker_Planck_collision}
    \mathcal{Q}[f] = \gradv\cdot \left(T\gradv f + (\bv-\bu) f\right)
    = T \gradv \cdot \left(\mathcal{M}[f] \gradv \left(\frac{f}{\mathcal{M}[f]}\right)\right),
\end{equation}
where $\mathcal{M}[f]$ is the Maxwellian distribution defined by
\begin{equation}
\label{eq_Maxwellian}
    \mathcal{M}[f] = \frac{n(t,\bx)}{(2\pi T(t,\bx))^{3/2}} \exp\left(-\frac{\vert \bv - \bu(t,\bx) \vert^2}{2T(t,\bx)}\right),
\end{equation}
with the macroscopic quantities---density, bulk velocity, and temperature---given by
\begin{equation}
\label{eq_macroscopic_quantities}
    n = \int_{\mathbb{R}^3} f \dd \bv, \quad
    \bu = \frac{1}{n} \int_{\mathbb{R}^3} f\bv \dd v, \quad
    T = \frac{1}{3n} \int_{\mathbb{R}^3} f\vert \bv - \bu \vert^2 \dd \bv.
\end{equation}

For the rest of this paper, we assume that the solution to Eq.~\eqref{eq_kinetic_equation} is homogeneous in the second and third spatial dimensions for  simplicity. Hence, $f=f(t,x,\bv)$ with $x\in \Omega_{x}\subset \mathbb{R}$ and $\bv = (v^{(1)},v^{(2)},v^{(3)})\in \mathbb{R}^3$. The proposed method, however, does not rely on this assumption and can be extended directly to multiple spatial dimensions.

Our method is based on the projector-splitting integrator introduced in \cite{lubich2015time}, implemented using a discretize-then-project (DtP) approach \cite{wang2026dynamical}. For  \eqref{eq_kinetic_equation}, this means that we first discretize both the physical and velocity spaces to obtain a semi-discretized ODE system, and then project it onto the tangent space of a fixed-rank velocity-only TT manifold. We subsequently apply the projector-splitting integrator to solve the resulting projected system, employing appropriate explicit or implicit time discretization in each sub-projection step. In \cite{wang2026dynamical}, all substeps are solved using explicit time-stepping schemes; consequently, the method is subject to severe stability constraints. This limitation is particularly pronounced in the presence of the diffusive Fokker–Planck collision operator and for large collision strength $\eta$, as an explicit scheme would typically require a time step of $\mathcal{O}(\Delta v^2/\eta)$ for stability, where $\Delta v$ is the velocity mesh size. To alleviate this constraint, the main contribution of this work is the development of a simple and efficient dynamical TT algorithm that treats the Fokker--Planck operator implicitly. To this end, we first introduce the spatial discretization and the TT representation of the Fokker--Planck operator.

\subsection{Spatial discretization and TT representation of the Fokker--Planck operator}

We first introduce the grid points in phase space. Without loss of generality, we assume the physical domain $\Omega_x=[0,L_x]$ with grid points $x_j=(j-1/2)\Delta x$, $j=1,\dots,N_x$, where $\Delta x=L_x/N_x$. For the velocity space, we truncate it to a sufficiently large domain $[v_{\min},v_{\max}]^3$ and  choose grid points $v_k=v_{\min}+(k-1/2)\Delta v$, $k=1,\dots,N_v$, with $\Delta v=(v_{\max}-v_{\min})/N_v$ in each velocity dimension. 

With these notations, the discrete value of $f$ on the grid point $(x_j,v_{k_1},v_{k_2},v_{k_3})$ at time $t_n=n\dt$ is given by
\begin{equation}
    f_{k_1k_2k_3}^{j,n} \approx f(t_n,x_j,v_{k_1},v_{k_2},v_{k_3}).
\end{equation}
At each fixed time $t_n$ and spatial grid point $x_j$, we represent the three-way tensor $\bf:=\{f_{k_1k_2k_3}^{j,n}\}\in \mathbb{R}^{N_v\times N_v\times N_v}$ in TT format (suppressing the indices $j$ and $n$ for brevity):
\begin{equation}
\label{eq_general_tensor_train}
    f_{k_1k_2k_3} = \sum_{\alpha_1=1}^{r_1}
    \sum_{\alpha_2=1}^{r_2}
    f_{k_1\alpha_1}^{(1)}
    f_{\alpha_1 k_2\alpha_2}^{(2)}
    f_{\alpha_2 k_3}^{(3)},
\end{equation}
where $f^{(1)}\in \mathbb{R}^{N_v\times r_1}$, $f^{(2)}\in \mathbb{R}^{r_1\times N_v\times r_2}$, $f^{(3)}\in \mathbb{R}^{r_2\times N_v}$ are the tensor cores of $\bf$, and $(r_1,r_2)$ is its TT-rank. For notational convenience, we use 
\begin{equation}
    \bf = \llbracket f^{(1)}, f^{(2)}, f^{(3)} \rrbracket
\end{equation}
to denote the TT representation of $\bf$ with tensor cores $f^{(1)},f^{(2)},f^{(3)}$, and we also use the following tensor diagram:
\begin{equation}
\label{diag_general_tensor_train}
\begin{tikzpicture}[baseline=0]
\draw (0,0) -- (1.6,0);
\draw (1.6,0) -- (3.2,0);
\draw (0,0) -- (0,-0.8);
\draw (1.6,0) -- (1.6,-0.8);
\draw (3.2,0) -- (3.2,-0.8);
\draw[fill=white,line width = 1.2] (0,0) circle (0.2);
\draw[fill=white,line width = 1.2] (1.6,0) circle (0.2);
\draw[fill=white,line width = 1.2] (3.2,0) circle (0.2);
\node at (0,0.5) {$f^{(1)}$};
\node at (1.6,0.5) {$f^{(2)}$};
\node at (3.2,0.5) {$f^{(3)}$};
\node at (0.8,0.2) {$\alpha_1$};
\node at (2.4,0.2) {$\alpha_2$};
\node at (0,-1) {$k_1$};
\node at (1.6,-1) {$k_2$};
\node at (3.2,-1) {$k_3$};
\end{tikzpicture}
\end{equation}
In this diagram, $f^{(1)}$, $f^{(2)}$, and $f^{(3)}$ have two, three, and two legs, respectively, corresponding to the number of indices of each tensor core in Eq.~\eqref{eq_general_tensor_train}.
The $\alpha_1$ leg of $f^{(1)}$ is connected to the $\alpha_1$ leg of $f^{(2)}$, representing a summation over this index; similarly, the $\alpha_2$ leg of $f^{(2)}$ is connected to that of $f^{(3)}$.
In this way, \eqref{diag_general_tensor_train} provides a visual representation of \eqref{eq_general_tensor_train} without explicit summation notation.

We are now ready to discuss the representation of the discrete Fokker--Planck operator in TT format.
We will show that, unlike a general TT operator, the Fokker--Planck  operator \eqref{eq_Fokker_Planck_collision} admits a simple structure under a certain finite difference discretization when the Maxwellian $\mathcal{M}$ is given.
This structure will be explored in this subsection and leveraged in the subsequent design of the time integrator.

First, we observe that the three velocity components in the Maxwellian distribution \eqref{eq_Maxwellian} are completely separable:
\begin{equation}
    \mathcal{M}[f] = \frac{n}{(2\pi T)^{3/2}}m^{(1)}(v^{(1)})
    m^{(2)}(v^{(2)})
    m^{(3)}(v^{(3)}),
\end{equation}
where
\begin{equation}
    m^{(l)}(v^{(l)})
    = \exp\left(-\frac{(v^{(l)}-u^{(l)})^2}{2T}\right), \quad l=1,2,3.
\end{equation}
On the grid point $(v_{k_1},v_{k_2},v_{k_3})$, the function value is therefore given by
\begin{equation}
    \mathcal{M}_{k_1k_2k_3} = \frac{n}{(2\pi T)^{3/2}} m^{(1)}_{k_1} m^{(2)}_{k_2}
    m^{(3)}_{k_3},
\end{equation}
which is readily a tensor train $\bM:=\{\mathcal{M}_{k_1k_2k_3}\}\in\mathbb{R}^{N_v\times N_v\times N_v}$ with TT-rank $(1,1)$. Note that here and in the following we again suppress the dependence on $t$ and $x$ (i.e., the indices $j$ and $n$) since the Fokker--Planck operator acts only in velocity space.

Given the Maxwellian $\bM$, we then split the Fokker--Planck operator (using the second form in \eqref{eq_Fokker_Planck_collision}) into three directions:
\begin{equation}
\label{eq_Fokker_Planck}
\bQ_{\bM}[\bf] = \bQ^{(1)}_{\bM}[\bf]
+ \bQ^{(2)}_{\bM}[\bf] + \bQ^{(3)}_{\bM}[\bf],
\end{equation}
where the first direction is discretized using a second-order central finite difference:
\begin{equation}
\label{eq_Fokker_Planck1}
    \left(\bQ^{(1)}_{\bM}[\bf]\right)_{k_1k_2k_3} = \frac{T}{\dv}\left(\mathfrak{F}_{k_1+\frac{1}{2},k_2k_3} 
        -  \mathfrak{F}_{k_1-\frac{1}{2},k_2k_3}\right),
\end{equation}
with fluxes at the interior grid points $k_1=1,\cdots,N_v-1$ given by
\begin{equation}
\label{eq_Fokker_Planck_fluxes}
\begin{aligned}
    \mathfrak{F}_{k_1+\frac{1}{2},k_2k_3}
=\frac{\mathcal{M}_{k_1k_2k_3}+\mathcal{M}_{k_1+1,k_2k_3}}{2\dv}
    \left(
        \frac{f_{k_1+1,k_2k_3}}{\mathcal{M}_{k_1+1,k_2k_3}}
        - \frac{f_{k_1k_2k_3}}{\mathcal{M}_{k_1k_2k_3}}
    \right),
\end{aligned}
\end{equation}
and zero at the boundaries $\mathfrak{F}_{\frac{1}{2},k_2k_3}=\mathfrak{F}_{N_v+\frac{1}{2},k_2k_3}=0$.
To enforce the boundary conditions, we may simply set $\mathcal{M}_{0,k_2k_3} = -\mathcal{M}_{1,k_2k_3}$ and $\mathcal{M}_{N_{v+1},k_2k_3}=-\mathcal{M}_{N_v,k_2k_3}$, so that \eqref{eq_Fokker_Planck_fluxes} applies to both interior and boundary points.
The discretization in the other two velocity directions are analogous.

Assuming $\bf$ is represented as in \eqref{eq_general_tensor_train}, a further simplification of \eqref{eq_Fokker_Planck1} yields
\begin{equation}
\label{eq_first_velocity_direction}
\begin{aligned}
    \left(\bQ^{(1)}_{\bM}[\bf] \right)_{k_1k_2k_3}
    &= \frac{T}{2\dv^2}
    \Bigg(
    \left(1+\frac{\mathcal{M}_{k_1k_2k_3}}{\mathcal{M}_{k_1-1,k_2k_3}}\right) f_{k_1-1,k_2k_3} 
    +\left(1+\frac{\mathcal{M}_{k_1k_2k_3}}{\mathcal{M}_{k_1+1,k_2k_3}}\right) f_{k_1+1,k_2k_3}
    \\
    &\qquad -\left(2+\frac{\mathcal{M}_{k_1-1,k_2k_3}+\mathcal{M}_{k_1+1,k_2k_3}}{\mathcal{M}_{k_1k_2k_3}}\right) f_{k_1,k_2k_3} 
    \Bigg) \\
    &= \frac{T}{2\dv^2} 
    \Bigg(
        \left(1+\frac{m^{(1)}_{k_1}}{m^{(1)}_{k_1-1}}\right)f_{k_1-1,k_2k_3}
        + \left(1+\frac{m^{(1)}_{k_1}}{m^{(1)}_{k_1+1}}\right)f_{k_1+1,k_2k_3} \\
        &\qquad - \left(2+\frac{m^{(1)}_{k_1-1}+m^{(1)}_{k_1+1}}{m^{(1)}_{k_1}}\right)f_{k_1k_2k_3}
    \Bigg) \\
    &= \sum_{\alpha_1=1}^{r_1}
    \sum_{\alpha_2=1}^{r_2}
    \left(a_{k_1}^{(1)} f_{k_1-1,\alpha_1}^{(1)}
    + b_{k_1}^{(1)} f_{k_1+1,\alpha_1}^{(1)}
    + c_{k_1}^{(1)} f_{k_1\alpha_1}^{(1)}
    \right)f_{\alpha_1k_2\alpha_2}^{(2)}f_{\alpha_2k_3}^{(3)},
\end{aligned}
\end{equation}
with 
\begin{equation}
    a_{k_1}^{(1)} = \frac{T}{2\dv^2}\left(1+\frac{m_{k_1}^{(1)}}{m_{k_1-1}^{(1)}}\right),\quad
    b_{k_1}^{(1)} = \frac{T}{2\dv^2}\left(1+\frac{m_{k_1}^{(1)}}{m_{k_1+1}^{(1)}}\right), \quad
    c_{k_1}^{(1)} = -\frac{T}{2\dv^2}\left(2 + \frac{m_{k_1-1}^{(1)}+m_{k_1+1}^{(1)}}{m_{k_1}^{(1)}}\right).
\end{equation}
Note that $a_{1}^{(1)} = b_{N_v}^{(1)} = 0$ due to the zero-flux boundary condition. 
From \eqref{eq_first_velocity_direction}, we can see that $\bQ^{(1)}_{\bM}[\bf]$ remains in TT form, with only the first core modified, while the second and third cores remain unchanged. Specifically, it can be represented as
\begin{equation}
\label{eq_QM1f}
    \bQ^{(1)}_{\bM}[\bf] = \llbracket  f^{(1)}\times_1 J^{(1)}, f^{(2)},f^{(3)}\rrbracket,
\end{equation} %
where $J^{(1)}$ is a tridiagonal matrix of size $N_v\times N_v$ given by
\begin{equation}
\label{eq_tridiagonal_matrix}
J^{(1)} = {\renewcommand{\arraystretch}{1.4}
\begin{pmatrix}
c_1^{(1)} & a_2^{(1)} & & & \\
b_1^{(1)} & c_2^{(1)} & a_3^{(1)} & & \\
& b_2^{(1)} & \ddots & \ddots & \\
& & \ddots & \ddots & a_{N_v}^{(1)} \\
& & & b_{N_v-1}^{(1)} & c_{N_v}^{(1)}
\end{pmatrix}
}.
\end{equation}
In \eqref{eq_QM1f}, $\times_1$ denotes a contraction between a tensor and a matrix, whose general form is defined as\footnote{We note that this definition differs by a transpose from the tensor $n$-mode product defined in \cite{kolda2009tensor}. We choose it for the convenience of the following discussion of the tensor Sylvester equation.}
\begin{equation}
    (X\times_n Y)_{i_1i_2\cdots i_{n-1}ji_{n+1}\dots i_d}
    = \sum_{i_n=1}^{I_n} X_{i_1\cdots i_{n-1}i_n i_{n+1}\cdots i_d}
    Y_{i_nj},
\end{equation}
where $X\in \mathbb{R}^{I_1\times I_2\times\cdots\times I_d}$, $Y\in \mathbb{R}^{I_n\times J}$, and $X\times_nY\in \mathbb{R}^{I_1\times \cdots \times I_{n-1}\times J\times I_{n+1}\times\cdots \times I_d}$.
Diagrammatically, $\bQ^{(1)}_{\bM}[\bf]$ can be represented as
\begin{equation}
\begin{tikzpicture}[baseline=0]
\draw (0,0) -- (1.6,0);
\draw (1.6,0) -- (3.2,0);
\draw (0,0) -- (0,-1.6);
\draw (1.6,0) -- (1.6,-0.8);
\draw (3.2,0) -- (3.2,-0.8);
\draw[fill=white,line width = 1.2] (0,0) circle (0.2);
\draw[fill=white,line width = 1.2] (1.6,0) circle (0.2);
\draw[fill=white,line width = 1.2] (3.2,0) circle (0.2);
\node at (0,0.5) {$f^{(1)}$};
\node at (1.6,0.5) {$f^{(2)}$};
\node at (3.2,0.5) {$f^{(3)}$};
\draw [fill=black] (0,-0.6) -- (-0.2,-0.8) -- (0,-1) -- (0.2,-0.8);
\node at (-0.6,-0.8) {$J^{(1)}$};
\end{tikzpicture}
\end{equation}
where
the diamond-shaped node denotes the tridiagonal matrix $J^{(1)}$.
The row index corresponds to the top leg and the column index corresponds to the bottom leg.
In the following, we will always use the diamond shape to represent a tridiagonal matrix in tensor diagrams.

Similar to the first direction, we define the tridiagonal matrices $J^{(2)}$ and $J^{(3)}$.
The second and third velocity directions can therefore be written as
\begin{equation}
\mathbf{Q}^{(2)}_{\bM}[\bf] = \llbracket f^{(1)}, f^{(2)}\times_2 J^{(2)}, f^{(3)} \rrbracket,\qquad
\mathbf{Q}^{(3)}_{\bM}[\bf] = \llbracket f^{(1)}, f^{(2)},  f^{(3)} \times_2 J^{(3)} \rrbracket,
\end{equation}
which can also be represented using tensor diagrams.
Combining all three directions, we obtain the diagrammatic  representation of the Fokker--Planck operator \eqref{eq_Fokker_Planck}:
\begin{equation}
\bQ_\bM[\bf] = 
\begin{tikzpicture}[baseline=0]
\draw (0,0) -- (1.6,0);
\draw (1.6,0) -- (3.2,0);
\draw (0,0) -- (0,-1.6);
\draw (1.6,0) -- (1.6,-0.8);
\draw (3.2,0) -- (3.2,-0.8);
\draw[fill=white,line width = 1.2] (0,0) circle (0.2);
\draw[fill=white,line width = 1.2] (1.6,0) circle (0.2);
\draw[fill=white,line width = 1.2] (3.2,0) circle (0.2);
\node at (0,0.5) {$f^{(1)}$};
\node at (1.6,0.5) {$f^{(2)}$};
\node at (3.2,0.5) {$f^{(3)}$};
\draw [fill=black] (0,-0.6) -- (-0.2,-0.8) -- (0,-1) -- (0.2,-0.8);
\node at (0.6,-0.8) {$J^{(1)}$};
\end{tikzpicture}
+
\begin{tikzpicture}[baseline=0]
\draw (0,0) -- (1.6,0);
\draw (1.6,0) -- (3.2,0);
\draw (0,0) -- (0,-0.8);
\draw (1.6,0) -- (1.6,-1.6);
\draw (3.2,0) -- (3.2,-0.8);
\draw[fill=white,line width = 1.2] (0,0) circle (0.2);
\draw[fill=white,line width = 1.2] (1.6,0) circle (0.2);
\draw[fill=white,line width = 1.2] (3.2,0) circle (0.2);
\node at (0,0.5) {$f^{(1)}$};
\node at (1.6,0.5) {$f^{(2)}$};
\node at (3.2,0.5) {$f^{(3)}$};
\draw [fill=black] (1.6,-0.6) -- (1.4,-0.8) -- (1.6,-1) -- (1.8,-0.8);
\node at (2.2,-0.8) {$J^{(2)}$};
\end{tikzpicture}
+
\begin{tikzpicture}[baseline=0]
\draw (0,0) -- (1.6,0);
\draw (1.6,0) -- (3.2,0);
\draw (0,0) -- (0,-0.8);
\draw (1.6,0) -- (1.6,-0.8);
\draw (3.2,0) -- (3.2,-1.6);
\draw[fill=white,line width = 1.2] (0,0) circle (0.2);
\draw[fill=white,line width = 1.2] (1.6,0) circle (0.2);
\draw[fill=white,line width = 1.2] (3.2,0) circle (0.2);
\node at (0,0.5) {$f^{(1)}$};
\node at (1.6,0.5) {$f^{(2)}$};
\node at (3.2,0.5) {$f^{(3)}$};
\draw [fill=black] (3.2,-0.6) -- (3,-0.8) -- (3.2,-1) -- (3.4,-0.8);
\node at (2.6,-0.8) {$J^{(3)}$};
\end{tikzpicture}
\end{equation}
This expression indicates that, for a given Maxwellian $\mathcal{M}$, the Fokker--Planck collision operator $\mathcal{Q}[f]$, in the discrete sense, can be represented as the sum of three terms, each modifying only one tensor core of $f$ through multiplication by a tridiagonal matrix.

\subsection{Time integration of tensor trains}

We now describe the time integration of Eq.~\eqref{eq_kinetic_equation} using a dynamical TT approximation. 
Since the Fokker--Planck operator \eqref{eq_Fokker_Planck_collision} depends on the Maxwellian $\mathcal{M}[f]$, it is nonlinear. To enable the implicit treatment, we first compute the Maxwellian $\bM^{j,n+1}$ at grid point $x_j$ and time $t_{n+1}$ before updating $\bf^{j,n+1}$. This is achieved using the well-known discrete moment update approach (see, e.g., \cite{hu2017asymptotic}). We begin by discretizing \eqref{eq_kinetic_equation} with a first-order IMEX scheme:
\begin{equation}
\label{eq:IMEX1}
    \frac{\bf^{j,n+1}-\bf^{j,n}}{\dt}
    + \left(v^{(1)} \bD_{x}^\up \bf\right)^{j,n}
    + \left(F^{(1)} \bD_{v^{(1)}}^\up \bf\right)^{j,n}
    = \eta \bQ_{\bM^{j,n+1}}[\bf^{j,n+1}],
\end{equation}
where $\left(v^{(1)} \bD_{x}^\up \bf\right)^{j,n}$ and $\left(F^{(1)} \bD_{v^{(1)}}^\up \bf\right)^{j,n}$ denote upwind discretizations of the transport terms. We then take the discrete moments $\sum_{k_1,k_2,k_3=1}^{N_v} \cdot (1,\bv_k,|\bv_k|^2)^\top\Delta v^3$ of \eqref{eq:IMEX1}. Using the conservation property of the Fokker--Planck operator, we obtain
\begin{equation}
\label{eq_moments}
\frac{\bU^{j,n+1}-\bU^{j,n}}{\dt}
+ \sum_{k_1,k_2,k_3=1}^{N_v} 
\left(\left(v^{(1)}\bD_{x}^\up \bf\right)^{j,n}
+ \left(F^{(1)} \bD_{v^{(1)}}^\up \bf\right)^{j,n}\right)
(1,\bv_k,|\bv_k|^2)^\top \dv^3 = 0,
\end{equation}
where
\begin{equation}
\bU := \sum_{k_1,k_2,k_3=1}^{N_v} \bf (1,\bv_k,|\bv_k|^2)^\top \dv^3
= (n, n\bu, n|\bu|^2 + 3nT)^\top.
\end{equation}
In this way, the macroscopic quantities $n^{j,n+1},\bu^{j,n+1}$ and $T^{j,n+1}$ can be computed without explicitly knowing $\bf^{j,n+1}$. The Maxwellian $\bM^{j,n+1}$ can therefore be constructed directly from these macroscopic quantities.
With this approach, the Maxwellian  $\bM^{j,n+1}$ is treated as known at the beginning of each time step.

To solve $\bf^{j,n+1}$ in TT format, we adopt the first-order projector-splitting integrator \cite{lubich2015time}. This requires a five-substep procedure to update the tensor cores of $\bf^{j,n+1}$ sequentially from left to right over each time step, where substep 1, 3, and 5 are forward projection steps, and substeps 2 and 4 are backward projection steps. We briefly describe this method below. Our presentation largely follows that in \cite{wang2026dynamical}. However, unlike \cite{wang2026dynamical}, where explicit schemes are used for all substeps, we apply implicit or IMEX schemes to substeps 1, 3, and 5. 

To facilitate the exposition, we introduce the following five canonical TT forms:
\begin{gather}
\begin{tikzpicture}[baseline=0]
\draw (0,0) -- (1.6,0);
\draw (1.6,0) -- (3.2,0);
\draw (0,0) -- (0,-0.8);
\draw (1.6,0) -- (1.6,-0.8);
\draw (3.2,0) -- (3.2,-0.8);
\draw[fill=white,line width=1.2] (0,0) circle (0.2);
\draw[fill=white,line width=1.2] (1.4,0) -- (1.7,-0.2) -- (1.7,0.2) -- (1.4,0);
\draw[fill=white,line width=1.2] (3.0,0) -- (3.3,-0.2) -- (3.3,0.2) -- (3.0,0);
\node at (0,0.5) {$C^{(1)}$};
\node at (1.6,0.5) {$Q^{(2)}$};
\node at (3.2,0.5) {$Q^{(3)}$};
\node at (0.8,-0.25) {$\alpha_1$};
\node at (2.4,-0.25) {$\alpha_2$};
\node at (0,-1) {$k_1$};
\node at (1.6,-1) {$k_2$};
\node at (3.2,-1) {$k_3$};
\node at (1.6,-1.6) {\rombracket{1}};
\end{tikzpicture}
\begin{tikzpicture}[baseline=0]
\draw (0,0) -- (1.6,0);
\draw (1.6,0) -- (3.2,0);
\draw (0,0) -- (0,-0.8);
\draw (1.6,0) -- (1.6,-0.8);
\draw (3.2,0) -- (3.2,-0.8);
\draw[fill=white,line width=1.2] (0.2,0) -- (-0.1,-0.2) -- (-0.1,0.2) -- (0.2,0);
\draw[fill=white,line width=1.2] (0.8,0) circle (0.2);
\draw[fill=white,line width=1.2] (1.4,0) -- (1.7,-0.2) -- (1.7,0.2) -- (1.4,0);
\draw[fill=white,line width=1.2] (3.0,0) -- (3.3,-0.2) -- (3.3,0.2) -- (3.0,0);
\node at (-0.2,0.5) {$P^{(1)}$};
\node at (0.8,0.5) {$S^{(1)}$};
\node at (1.8,0.5) {$Q^{(2)}$};
\node at (3.2,0.5) {$Q^{(3)}$};
\node at (0.4,-0.25) {$\alpha_1$};
\node at (1.2,-0.25) {$\alpha_1'$};
\node at (2.4,-0.25) {$\alpha_2$};
\node at (0,-1) {$k_1$};
\node at (1.6,-1) {$k_2$};
\node at (3.2,-1) {$k_3$};
\node at (1.6,-1.6) {\rombracket{2}};
\end{tikzpicture}\nonumber\\
\begin{tikzpicture}[baseline=0]
\draw (0,0) -- (1.6,0);
\draw (1.6,0) -- (3.2,0);
\draw (0,0) -- (0,-0.8);
\draw (1.6,0) -- (1.6,-0.8);
\draw (3.2,0) -- (3.2,-0.8);
\draw[fill=white,line width=1.2] (0.2,0) -- (-0.1,-0.2) -- (-0.1,0.2) -- (0.2,0);
\draw[fill=white,line width=1.2] (1.6,0) circle (0.2);
\draw[fill=white,line width=1.2] (3.0,0) -- (3.3,-0.2) -- (3.3,0.2) -- (3.0,0);
\node at (0,0.5) {$P^{(1)}$};
\node at (1.6,0.5) {$C^{(2)}$};
\node at (3.2,0.5) {$Q^{(3)}$};
\node at (0.8,-0.25) {$\alpha_1$};
\node at (2.4,-0.25) {$\alpha_2$};
\node at (0,-1) {$k_1$};
\node at (1.6,-1) {$k_2$};
\node at (3.2,-1) {$k_3$};
\node at (1.6,-1.6) {\rombracket{3}};
\end{tikzpicture}
\begin{tikzpicture}[baseline=0]
\draw (0,0) -- (1.6,0);
\draw (1.6,0) -- (3.2,0);
\draw (0,0) -- (0,-0.8);
\draw (1.6,0) -- (1.6,-0.8);
\draw (3.2,0) -- (3.2,-0.8);
\draw[fill=white,line width=1.2] (0.2,0) -- (-0.1,-0.2) -- (-0.1,0.2) -- (0.2,0);
\draw[fill=white,line width=1.2] (1.8,0) -- (1.5,-0.2) -- (1.5,0.2) -- (1.8,0);
\draw[fill=white,line width=1.2] (2.4,0) circle (0.2);
\draw[fill=white,line width=1.2] (3.0,0) -- (3.3,-0.2) -- (3.3,0.2) -- (3.0,0);
\node at (0,0.5) {$P^{(1)}$};
\node at (1.4,0.5) {$P^{(2)}$};
\node at (2.4,0.5) {$S^{(2)}$};
\node at (3.4,0.5) {$Q^{(3)}$};
\node at (0.8,-0.25) {$\alpha_1$};
\node at (2.0,-0.25) {$\alpha_2'$};
\node at (2.8,-0.25) {$\alpha_2$};
\node at (0,-1) {$k_1$};
\node at (1.6,-1) {$k_2$};
\node at (3.2,-1) {$k_3$};
\node at (1.6,-1.6) {\rombracket{4}};
\end{tikzpicture}
\begin{tikzpicture}[baseline=0]
\draw (0,0) -- (1.6,0);
\draw (1.6,0) -- (3.2,0);
\draw (0,0) -- (0,-0.8);
\draw (1.6,0) -- (1.6,-0.8);
\draw (3.2,0) -- (3.2,-0.8);
\draw[fill=white,line width=1.2] (0.2,0) -- (-0.1,-0.2) -- (-0.1,0.2) -- (0.2,0);
\draw[fill=white,line width=1.2] (1.8,0) -- (1.5,-0.2) -- (1.5,0.2) -- (1.8,0);
\draw[fill=white,line width=1.2] (3.2,0) circle (0.2);
\node at (0,0.5) {$P^{(1)}$};
\node at (1.6,0.5) {$P^{(2)}$};
\node at (3.2,0.5) {$C^{(3)}$};
\node at (0.8,-0.25) {$\alpha_1$};
\node at (2.4,-0.25) {$\alpha_2$};
\node at (0,-1) {$k_1$};
\node at (1.6,-1) {$k_2$};
\node at (3.2,-1) {$k_3$};
\node at (1.6,-1.6) {\rombracket{5}};
\end{tikzpicture}
\label{eq_TT_decomposition_forms1}
\end{gather}
Here, the circles denote general tensor cores, while the triangles represent cores satisfying orthonormality conditions, with the orientation indicating the  direction of orthonormalization. Specifically, we have 
\begin{equation}
\label{eq_orthogonal_properties}
\begin{aligned}
    &\sum_{k_1=1}^{N_1} 
    P_{k_1\alpha_1}^{(1)}
    P_{k_1\alpha_1'}^{(1)} 
    = \delta_{\alpha_1\alpha_1'}, 
    &\sum_{k_2=1}^{N_2} 
    \sum_{\alpha_1=1}^{r_1}
    P_{\alpha_1 k_2 \alpha_2}^{(2)}
    P_{\alpha_1 k_2 \alpha_2'}^{(2)}
    = \delta_{\alpha_2\alpha_2'}, \\
    &\sum_{k_3=1}^{N_3} Q_{\alpha_2k_3}^{(3)}Q_{\alpha_2'k_3}^{(3)}=\delta_{\alpha_2\alpha_2'},
    &\sum_{k_2=1}^{N_2} \sum_{\alpha_2'=1}^{r_2}
         Q_{\alpha_1k_2\alpha_2'}^{(2)}
         Q_{\alpha_1'k_2\alpha_2'}^{(2)} = \delta_{\alpha_1\alpha_1'}.
\end{aligned}
\end{equation}
These five forms can also be written as $\llbracket C^{(1)}, Q^{(2)},Q^{(3)}\rrbracket$,
$\llbracket P^{(1)}\times_2 S^{(1)}, Q^{(2)},Q^{(3)}\rrbracket$,
$\llbracket P^{(1)}, C^{(2)},Q^{(3)}\rrbracket$,
$\llbracket P^{(1)}, P^{(2)},S^{(2)}\times_2 Q^{(3)}\rrbracket$, and 
$\llbracket P^{(1)}, P^{(2)},C^{(3)}\rrbracket$, respectively.

We use the following IMEX scheme (or implicit scheme in the absense of the transport terms) in the forward projection steps:
\begin{equation}
\label{eq_forward_scheme}
    \frac{\tilde{\bf}^{j}-\bf^{j}}{\dt}
    + \left(v^{(1)} \bD_{x}^\up \bf\right)^{j}
    + \left(F^{(1)} \bD_{v^{(1)}}^\up \bf\right)^{j}
    = \eta \bQ_{\bM^{j,n+1}}[\tilde{\bf}^{j}],
\end{equation}
where $\bf^{j}$ is known, while $\tilde{\bf}^{j}$ is the unknown. This scheme is equivalent to
\begin{equation}
\label{eq_forward_scheme_transformed}
(\bI-\dt\eta \bQ_{\bM^{j,n+1}})\tilde{\bf}^{j}
= \bf^{j} - \dt \left(
\left(v^{(1)} \bD_{x}^\up \bf\right)^{j}
    + \left(F^{(1)}\bD_{v^{(1)}}^\up \bf\right)^{j}
\right).
\end{equation}
We use the following explicit scheme in the backward projection steps:
\begin{equation}
\label{eq_backward_scheme}
    \frac{\tilde{\bf}^{j}-{\bf}^{j}}{-\dt}
    + \left(v^{(1)} \bD_{x}^\up \bf\right)^{j}
    + \left(F^{(1)} \bD_{v^{(1)}}^\up \bf\right)^{j}
    = \eta \bQ_{\bM^{j,n+1}}[{\bf}^{j}],
\end{equation}
or equivalently,
\begin{equation}
\label{eq_backward_scheme_transformed}
    \tilde{\bf}^j = \bf^j + \dt \left(\left(v^{(1)} \bD_{x}^\up \bf\right)^{j}
    +\left(F^{(1)} \bD_{v^{(1)}}^\up \bf\right)^{j}\right)
    -\dt\eta \bQ_{\bM^{j,n+1}}[{\bf}^{j}].
\end{equation}

At time $t_n$ and each spatial point $x_j$, we assume  that $\bf^{j,n}$ is written in the form (I) above, denoted by $\bf_\rom{1}^{j}=\llbracket C^{(1)}, Q^{(2)},Q^{(3)}\rrbracket$. A first-order projector-splitting method then proceeds as follows: 
\begin{itemize}
    \item \textbf{Substep 1.} Taking $\bf^j=\bf_\rom{1}^{j}$, project \eqref{eq_forward_scheme_transformed} onto the cores $Q^{(2)}$ and $Q^{(3)}$ to obtain $\tilde{C}^{(1)}$.
    Update the solution to $\tilde{\bf}_\rom{1}^j=\llbracket \tilde{C}^{(1)}, Q^{(2)},Q^{(3)}\rrbracket$. Then convert $\tilde{\bf}_\rom{1}^{j}$ to form $\rombracket{2}$ by a QR decomposition on $\tilde{C}^{(1)}$, and denote the resulting solution by $\bf_\rom{2}^{j}=\llbracket P^{(1)}\times_2 S^{(1)}, Q^{(2)},Q^{(3)}\rrbracket$.
    \item \textbf{Substep 2.} Taking $\bf^j=\tilde{\bf}_\rom{1}^j$, project \eqref{eq_backward_scheme_transformed} onto the cores $P^{(1)}$, $Q^{(2)}$, and $Q^{(3)}$ to obtain $\tilde{S}^{(1)}$. Update the solution to $\tilde{\bf}_\rom{2}^j=\llbracket P^{(1)}\times_2 \tilde{S}^{(1)}, Q^{(2)},Q^{(3)}\rrbracket$. Then convert $\tilde{\bf}_\rom{2}^{j}$ to form $\rombracket{3}$ by multiplying $\tilde{S}^{(1)}$ and $Q^{(2)}$, and denote the resulting solution by  $\bf_\rom{3}^{j}=\llbracket P^{(1)},C^{(2)},Q^{(3)}\rrbracket$.
    \item \textbf{Substep 3.} Taking $\bf^j=\bf_\rom{3}^{j}$, project \eqref{eq_forward_scheme_transformed} onto the cores $P^{(1)}$ and $Q^{(3)}$ to obtain $\tilde{C}^{(2)}$. Update the solution to $\tilde{\bf}_\rom{3}^j=\llbracket P^{(1)},\tilde{C}^{(2)},Q^{(3)}\rrbracket$. Then convert $\tilde{\bf}_\rom{3}^{j}$ to form $\rombracket{4}$ by a QR decomposition on $\tilde{C}^{(2)}$, and denote the resulting solution by $\bf_\rom{4}^{j}=\llbracket P^{(1)}, P^{(2)},S^{(2)}\times_2 Q^{(3)}\rrbracket$.
    \item \textbf{Substep 4.} Taking $\bf^j=\tilde{\bf}_\rom{3}^{j}$, project \eqref{eq_backward_scheme_transformed} onto the cores $P^{(1)}$, $P^{(2)}$, and $Q^{(3)}$ to obtain $\tilde{S}^{(2)}$. Update the solution to $\tilde{\bf}_\rom{4}^j=\llbracket P^{(1)}, P^{(2)},\tilde{S}^{(2)}\times_2 Q^{(3)}\rrbracket$. Then convert $\tilde{\bf}_\rom{4}^{j}$ to form $\rombracket{5}$ by multiplying $\tilde{S}^{(2)}$ and $Q^{(3)}$, and denote the resulting solution by  $\bf_\rom{5}^{j}=\llbracket P^{(1)},P^{(2)},C^{(3)}\rrbracket$.
    \item \textbf{Substep 5.} Taking $\bf^j=\bf_\rom{5}^{j}$, project \eqref{eq_forward_scheme_transformed} onto the cores $P^{(1)}$ and $P^{(2)}$ to obtain $\tilde{C}^{(3)}$. Update the solution to $\tilde{\bf}_\rom{5}^j=\llbracket P^{(1)},{P}^{(2)},\tilde{C}^{(3)}\rrbracket$. Then convert $\tilde{\bf}_\rom{5}^{j}$ to form $\rombracket{1}$ by a left-orthonormalization procedure, yielding the new solution $\bf^{j,n+1}$.
\end{itemize}

In the above algorithm, the projections in \textbf{substeps 2} and \textbf{4} are relatively straightforward, since the underlying scheme \eqref{eq_backward_scheme_transformed} is fully explicit (more details can be found in \cite{wang2026dynamical}). The projections in \textbf{substeps 1, 3}, and \textbf{5} are more involved, since the underlying scheme \eqref{eq_forward_scheme_transformed} is implicit. Addressing this issue is the focus of the next subsection.

\subsection{Implicit treatment of substeps 1, 3, and 5}

We start from the diagrammatical representation of Eq.~\eqref{eq_forward_scheme_transformed}:
\begin{equation}
\begin{aligned}
\label{diag_forward_scheme}
\begin{tikzpicture}[baseline=-18]
\draw (0,0) -- (1.2,0);
\draw (1.2,0) -- (2.4,0);
\draw (0,0) -- (0,-1.6);
\draw (1.2,0) -- (1.2,-0.8);
\draw (2.4,0) -- (2.4,-0.8);
\draw[fill=white,line width = 1.2] (0,0) circle (0.2);
\draw[fill=white,line width = 1.2] (1.2,0) circle (0.2);
\draw[fill=white,line width = 1.2] (2.4,0) circle (0.2);
\node at (0,0.5) {$\tilde{f}^{(1)}$};
\node at (1.2,0.5) {$\tilde{f}^{(2)}$};
\node at (2.4,0.5) {$\tilde{f}^{(3)}$};
\draw [fill=black] (0,-0.6) -- (-0.2,-0.8) -- (0,-1) -- (0.2,-0.8);
\node at (0.6,-0.8) {$T^{(1)}$};
\end{tikzpicture}
+
\begin{tikzpicture}[baseline=-18]
\draw (0,0) -- (1.2,0);
\draw (1.2,0) -- (2.4,0);
\draw (0,0) -- (0,-0.8);
\draw (1.2,0) -- (1.2,-1.6);
\draw (2.4,0) -- (2.4,-0.8);
\draw[fill=white,line width = 1.2] (0,0) circle (0.2);
\draw[fill=white,line width = 1.2] (1.2,0) circle (0.2);
\draw[fill=white,line width = 1.2] (2.4,0) circle (0.2);
\node at (0,0.5) {$\tilde{f}^{(1)}$};
\node at (1.2,0.5) {$\tilde{f}^{(2)}$};
\node at (2.4,0.5) {$\tilde{f}^{(3)}$};
\draw [fill=black] (1.2,-0.6) -- (1.,-0.8) -- (1.2,-1) -- (1.4,-0.8);
\node at (1.8,-0.8) {$T^{(2)}$};
\end{tikzpicture}
+
\begin{tikzpicture}[baseline=-18]
\draw (0,0) -- (1.2,0);
\draw (1.2,0) -- (2.4,0);
\draw (0,0) -- (0,-0.8);
\draw (1.2,0) -- (1.2,-0.8);
\draw (2.4,0) -- (2.4,-1.6);
\draw[fill=white,line width = 1.2] (0,0) circle (0.2);
\draw[fill=white,line width = 1.2] (1.2,0) circle (0.2);
\draw[fill=white,line width = 1.2] (2.4,0) circle (0.2);
\node at (0,0.5) {$\tilde{f}^{(1)}$};
\node at (1.2,0.5) {$\tilde{f}^{(2)}$};
\node at (2.4,0.5) {$\tilde{f}^{(3)}$};
\draw [fill=black] (2.4,-0.6) -- (2.2,-0.8) -- (2.4,-1) -- (2.6,-0.8);
\node at (1.8,-0.8) {$T^{(3)}$};
\end{tikzpicture}
= 
\begin{tikzpicture}[baseline=-18]
\draw (0,0) -- (1.2,0);
\draw (1.2,0) -- (2.4,0);
\draw (0,0) -- (0,-0.8);
\draw (1.2,0) -- (1.2,-0.8);
\draw (2.4,0) -- (2.4,-0.8);
\draw[fill=white,line width = 1.2] (0,0) circle (0.2);
\draw[fill=white,line width = 1.2] (1.2,0) circle (0.2);
\draw[fill=white,line width = 1.2] (2.4,0) circle (0.2);
\node at (0,0.5) {${K}^{(1)}$};
\node at (1.2,0.5) {${K}^{(2)}$};
\node at (2.4,0.5) {${K}^{(3)}$};
\end{tikzpicture}
\end{aligned}
\end{equation}
where $T^{(d)} := \frac{1}{3}I-\dt\eta J^{(d)}$ for $d=1,2,3$.
In this diagram, we use $\llbracket \tilde{f}^{(1)},\tilde{f}^{(2)},\tilde{f}^{(3)}\rrbracket$ to represent the unknown quantity $\tilde{\bf}^{j}$.
The tensor train $\bK=\llbracket K^{(1)}, K^{(2)},K^{(3)}\rrbracket$ represents the entire right-hand side of \eqref{eq_forward_scheme_transformed}.
In this formulation, we distribute the identity operator evenly among the three directions, which gives rise to the factor $\frac{1}{3}I$ in the definition of $T^{(d)}$. This choice is made purely for convenience; alternatively, the identity operator could be assigned entirely to a single velocity direction or distributed in other ways without affecting the following discussion.

As described in the previous subsection, in \textbf{substep 1}, the unknown tensor train $\tilde{\bf}^j$ is in the form \rombracket{1}, denoted by $\llbracket \tilde{C}^{(1)},Q^{(2)},Q^{(3)} \rrbracket$, and only the first core $\tilde{C}^{(1)}$ is updated. This update is performed through the projection of \eqref{eq_forward_scheme_transformed} onto the cores $Q^{(2)}$ and $Q^{(3)}$. 
Diagrammatically, this means the whole equation is contracted with $Q^{(2)}$ and $Q^{(3)}$:
\begin{equation}
\label{diag_substep1_projection}
\underbrace{
\begin{tikzpicture}[baseline=-26]
\draw[rounded corners=2pt, dashed, gray] (0.8,1) rectangle (2.8,-3);
\draw (0,0) -- (1.2,0);
\draw (1.2,0) -- (2.4,0);
\draw (0,0) -- (0,-1.6);
\draw (1.2,0) -- (1.2,-2);
\draw (2.4,0) -- (2.4,-2);
\draw[fill=white,line width = 1.2] (0,0) circle (0.2);
\draw[fill=white,line width=1.2] (1.0,0) -- (1.3,-0.2) -- (1.3,0.2) -- (1.0,0);
\draw[fill=white,line width=1.2] (2.2,0) -- (2.5,-0.2) -- (2.5,0.2) -- (2.2,0);
\node at (0,0.5) {$\tilde{C}^{(1)}$};
\node at (1.2,0.5) {${Q}^{(2)}$};
\node at (2.4,0.5) {${Q}^{(3)}$};
\draw [fill=black] (0,-0.8) -- (-0.2,-1) -- (0,-1.2) -- (0.2,-1);
\node at (0.6,-1) {$T^{(1)}$};
\draw (0.4,-2) -- (2.4,-2);
\draw[fill=white,line width=1.2] (1.0,-2) -- (1.3,-2.2) -- (1.3,-1.8) -- (1.0,-2);
\draw[fill=white,line width=1.2] (2.2,-2) -- (2.5,-2.2) -- (2.5,-1.8) -- (2.2,-2);
\node at (1.2,-2.5) {${Q}^{(2)}$};
\node at (2.4,-2.5) {${Q}^{(3)}$};
\end{tikzpicture}
}_{\tilde{C}^{(1)}\times_1 T^{(1)}}
+
\underbrace{
\begin{tikzpicture}[baseline=-26]
\draw[rounded corners=2pt, dashed, gray] (0.8,1) rectangle (2.8,-3);
\draw (0,0) -- (1.2,0);
\draw (1.2,0) -- (2.4,0);
\draw (0,0) -- (0,-1.6);
\draw (1.2,0) -- (1.2,-2);
\draw (2.4,0) -- (2.4,-2);
\draw[fill=white,line width = 1.2] (0,0) circle (0.2);
\draw[fill=white,line width=1.2] (1.0,0) -- (1.3,-0.2) -- (1.3,0.2) -- (1.0,0);
\draw[fill=white,line width=1.2] (2.2,0) -- (2.5,-0.2) -- (2.5,0.2) -- (2.2,0);
\node at (0,0.5) {$\tilde{C}^{(1)}$};
\node at (1.2,0.5) {${Q}^{(2)}$};
\node at (2.4,0.5) {${Q}^{(3)}$};
\draw [fill=black] (1.2,-0.8) -- (1,-1) -- (1.2,-1.2) -- (1.4,-1);
\node at (1.8,-1) {$T^{(2)}$};
\draw (0.4,-2) -- (2.4,-2);
\draw[fill=white,line width=1.2] (1.0,-2) -- (1.3,-2.2) -- (1.3,-1.8) -- (1.0,-2);
\draw[fill=white,line width=1.2] (2.2,-2) -- (2.5,-2.2) -- (2.5,-1.8) -- (2.2,-2);
\node at (1.2,-2.5) {${Q}^{(2)}$};
\node at (2.4,-2.5) {${Q}^{(3)}$};
\end{tikzpicture}
+
\begin{tikzpicture}[baseline=-26]
\draw[rounded corners=2pt, dashed, gray] (0.8,1) rectangle (2.8,-3);
\draw (0,0) -- (1.2,0);
\draw (1.2,0) -- (2.4,0);
\draw (0,0) -- (0,-1.6);
\draw (1.2,0) -- (1.2,-2);
\draw (2.4,0) -- (2.4,-2);
\draw[fill=white,line width = 1.2] (0,0) circle (0.2);
\draw[fill=white,line width=1.2] (1.0,0) -- (1.3,-0.2) -- (1.3,0.2) -- (1.0,0);
\draw[fill=white,line width=1.2] (2.2,0) -- (2.5,-0.2) -- (2.5,0.2) -- (2.2,0);
\node at (0,0.5) {$\tilde{C}^{(1)}$};
\node at (1.2,0.5) {${Q}^{(2)}$};
\node at (2.4,0.5) {${Q}^{(3)}$};
\draw [fill=black] (2.4,-0.8) -- (2.2,-1) -- (2.4,-1.2) -- (2.6,-1);
\node at (1.8,-1) {$T^{(3)}$};
\draw (0.4,-2) -- (2.4,-2);
\draw[fill=white,line width=1.2] (1.0,-2) -- (1.3,-2.2) -- (1.3,-1.8) -- (1.0,-2);
\draw[fill=white,line width=1.2] (2.2,-2) -- (2.5,-2.2) -- (2.5,-1.8) -- (2.2,-2);
\node at (1.2,-2.5) {${Q}^{(2)}$};
\node at (2.4,-2.5) {${Q}^{(3)}$};
\end{tikzpicture}
}_{\tilde{C}^{(1)}\times_2 H_\rom{1}}
=
\underbrace{
\begin{tikzpicture}[baseline=-15]
\draw (0,0) -- (1.2,0);
\draw (1.2,0) -- (2.4,0);
\draw (0.4,-1.2) -- (1.2,-1.2);
\draw (1.2,-1.2) -- (2.4,-1.2);
\draw (0,0) -- (0,-0.8);
\draw (1.2,0) -- (1.2,-1.2);
\draw (2.4,0) -- (2.4,-1.2);
\draw[fill=white,line width = 1.2] (0,0) circle (0.2);
\draw[fill=white,line width = 1.2] (1.2,0) circle (0.2);
\draw[fill=white,line width = 1.2] (2.4,0) circle (0.2);
\draw[fill=white,line width=1.2] (1.0,-1.2) -- (1.3,-1.4) -- (1.3,-1) -- (1.0,-1.2);
\draw[fill=white,line width=1.2] (2.2,-1.2) -- (2.5,-1.4) -- (2.5,-1) -- (2.2,-1.2);
\node at (0,0.5) {${K}^{(1)}$};
\node at (1.2,0.5) {${K}^{(2)}$};
\node at (2.4,0.5) {${K}^{(3)}$};
\node at (1.2,-1.7) {${Q}^{(2)}$};
\node at (2.4,-1.7) {${Q}^{(3)}$};
\end{tikzpicture}
}_{R_\rom{1}}
\end{equation}
Using the orthogonality property \eqref{eq_orthogonal_properties}, this yields a matrix Sylvester equation for $\tilde{C}^{(1)}$ of the form:
\begin{equation}
\label{eq_Sylvester_substep1}
    \tilde{C}^{(1)} \times_1 T^{(1)} + \tilde{C}^{(1)} \times_2 H_\rom{1} = R_\rom{1},
\end{equation}
where the unknown $\tilde{C}^{(1)}$ and the right-hand side $R_\rom{1}$ are matrices of size $N_v \times r_1$;
$T^{(1)}$ is a tridiagonal matrix of size $N_v\times N_v$, while $H_\rom{1}$ is a matrix of size $r_1\times r_1$.
In principle, any Sylvester solver can be used to solve Eq.~\eqref{eq_Sylvester_substep1} for $\tilde{C}^{(1)}$.
However, different solvers have different computational complexities, which we will discuss in the next section.

In \textbf{substep 3}, the unknown tensor train $\tilde{\bf}^j$ is in the form \rombracket{3}, denoted by $\llbracket P^{(1)},\tilde{C}^{(2)},Q^{(3)}\rrbracket$, and only the second core $\tilde{C}^{(2)}$ is updated. This update is performed through the projection of \eqref{eq_forward_scheme_transformed} onto the cores $P^{(1)}$ and $Q^{(3)}$. Diagrammatically, this means that the whole equation is contracted with $P^{(1)}$ and $Q^{(3)}$:
\begin{equation}
\label{diag_substep3_projection}
\underbrace{
\begin{tikzpicture}[baseline=-26]
\draw[rounded corners=2pt, dashed, gray] (-0.5,1) rectangle (0.5,-3);
\draw[rounded corners=2pt, dashed, gray] (1.9,1) rectangle (2.9,-3);
\draw (0,0) -- (1.2,0);
\draw (1.2,0) -- (2.4,0);
\draw (0,-2) -- (0.8,-2);
\draw (1.6,-2) -- (2.4,-2);
\draw (0,0) -- (0,-2);
\draw (1.2,0) -- (1.2,-1.6);
\draw (2.4,0) -- (2.4,-2);
\draw[fill=white,line width=1.2] (-0.1,0.2) -- (-0.1,-0.2) -- (0.2,0) -- (-0.1,0.2);
\draw[fill=white,line width = 1.2] (1.2,0) circle (0.2);
\draw[fill=white,line width=1.2] (2.2,0) -- (2.5,-0.2) -- (2.5,0.2) -- (2.2,0);
\draw[fill=white,line width=1.2] (-0.1,-1.8) -- (-0.1,-2.2) -- (0.2,-2.0) -- (-0.1,-1.8);
\draw[fill=white,line width=1.2] (2.2,-2) -- (2.5,-2.2) -- (2.5,-1.8) -- (2.2,-2);
\node at (0,0.5) {${P}^{(1)}$};
\node at (1.2,0.5) {$\tilde{C}^{(2)}$};
\node at (2.4,0.5) {${Q}^{(3)}$};
\node at (0,-2.5) {${P}^{(1)}$};
\node at (2.4,-2.5) {${Q}^{(3)}$};
\draw [fill=black] (0,-0.8) -- (-0.2,-1) -- (0,-1.2) -- (0.2,-1);
\node at (0.6,-1) {$T^{(1)}$};
\end{tikzpicture}
}_{\tilde{C}^{(2)}\times_1 G_\rom{3}}
+
\underbrace{
\begin{tikzpicture}[baseline=-26]
\draw[rounded corners=2pt, dashed, gray] (-0.5,1) rectangle (0.5,-3);
\draw[rounded corners=2pt, dashed, gray] (1.9,1) rectangle (2.9,-3);
\draw (0,0) -- (1.2,0);
\draw (1.2,0) -- (2.4,0);
\draw (0,-2) -- (0.8,-2);
\draw (1.6,-2) -- (2.4,-2);
\draw (0,0) -- (0,-2);
\draw (1.2,0) -- (1.2,-1.6);
\draw (2.4,0) -- (2.4,-2);
\draw[fill=white,line width=1.2] (-0.1,0.2) -- (-0.1,-0.2) -- (0.2,0) -- (-0.1,0.2);
\draw[fill=white,line width = 1.2] (1.2,0) circle (0.2);
\draw[fill=white,line width=1.2] (2.2,0) -- (2.5,-0.2) -- (2.5,0.2) -- (2.2,0);
\draw[fill=white,line width=1.2] (-0.1,-1.8) -- (-0.1,-2.2) -- (0.2,-2.0) -- (-0.1,-1.8);
\draw[fill=white,line width=1.2] (2.2,-2) -- (2.5,-2.2) -- (2.5,-1.8) -- (2.2,-2);
\node at (0,0.5) {${P}^{(1)}$};
\node at (1.2,0.5) {$\tilde{C}^{(2)}$};
\node at (2.4,0.5) {${Q}^{(3)}$};
\node at (0,-2.5) {${P}^{(1)}$};
\node at (2.4,-2.5) {${Q}^{(3)}$};
\draw [fill=black] (1.2,-0.8) -- (1,-1) -- (1.2,-1.2) -- (1.4,-1);
\node at (1.8,-1) {$T^{(2)}$};
\end{tikzpicture}
}_{\tilde{C}^{(2)}\times_2 T^{(2)}}
+
\underbrace{
\begin{tikzpicture}[baseline=-26]
\draw[rounded corners=2pt, dashed, gray] (-0.5,1) rectangle (0.5,-3);
\draw[rounded corners=2pt, dashed, gray] (1.9,1) rectangle (2.9,-3);
\draw (0,0) -- (1.2,0);
\draw (1.2,0) -- (2.4,0);
\draw (0,-2) -- (0.8,-2);
\draw (1.6,-2) -- (2.4,-2);
\draw (0,0) -- (0,-2);
\draw (1.2,0) -- (1.2,-1.6);
\draw (2.4,0) -- (2.4,-2);
\draw[fill=white,line width=1.2] (-0.1,0.2) -- (-0.1,-0.2) -- (0.2,0) -- (-0.1,0.2);
\draw[fill=white,line width = 1.2] (1.2,0) circle (0.2);
\draw[fill=white,line width=1.2] (2.2,0) -- (2.5,-0.2) -- (2.5,0.2) -- (2.2,0);
\draw[fill=white,line width=1.2] (-0.1,-1.8) -- (-0.1,-2.2) -- (0.2,-2.0) -- (-0.1,-1.8);
\draw[fill=white,line width=1.2] (2.2,-2) -- (2.5,-2.2) -- (2.5,-1.8) -- (2.2,-2);
\node at (0,0.5) {${P}^{(1)}$};
\node at (1.2,0.5) {$\tilde{C}^{(2)}$};
\node at (2.4,0.5) {${Q}^{(3)}$};
\node at (0,-2.5) {${P}^{(1)}$};
\node at (2.4,-2.5) {${Q}^{(3)}$};
\draw [fill=black] (2.4,-0.8) -- (2.2,-1) -- (2.4,-1.2) -- (2.6,-1);
\node at (1.8,-1) {$T^{(3)}$};
\end{tikzpicture}
}_{\tilde{C}^{(2)}\times_3 H_\rom{3}}
= 
\underbrace{
\begin{tikzpicture}[baseline=-14]
\draw (0,0) -- (1.2,0);
\draw (1.2,0) -- (2.4,0);
\draw (0,-1.2) -- (0.8,-1.2);
\draw (1.6,-1.2) -- (2.4,-1.2);
\draw (0,0) -- (0,-1.2);
\draw (1.2,0) -- (1.2,-0.8);
\draw (2.4,0) -- (2.4,-1.2);
\draw[fill=white,line width = 1.2] (0,0) circle (0.2);
\draw[fill=white,line width = 1.2] (1.2,0) circle (0.2);
\draw[fill=white,line width = 1.2] (2.4,0) circle (0.2);
\node at (0,0.5) {${K}^{(1)}$};
\node at (1.2,0.5) {${K}^{(2)}$};
\node at (2.4,0.5) {${K}^{(3)}$};
\draw[fill=white,line width=1.2] (0.2,-1.2) -- (-0.1,-1.4) -- (-0.1,-1) -- (0.2,-1.2);
\draw[fill=white,line width=1.2] (2.2,-1.2) -- (2.5,-1.4) -- (2.5,-1) -- (2.2,-1.2);
\node at (0,-1.8) {${P}^{(1)}$};
\node at (2.4,-1.8) {${Q}^{(3)}$};
\end{tikzpicture}
}_{R_\rom{3}}
\end{equation}
This leads to a tensor Sylvester equation for $\tilde{C}^{(2)}$:
\begin{equation}
\label{eq_Sylvester_substep3}
    \tilde{C}^{(2)}\times_1 G_\rom{3} + 
    \tilde{C}^{(2)}\times_2 T^{(2)} + 
    \tilde{C}^{(2)}\times_3 H_\rom{3} = R_\rom{3}.
\end{equation}
In \eqref{eq_Sylvester_substep3}, the unknown $\tilde{C}^{(2)}$ and the right-hand side $R_\rom{3}$ are tensors of size $r_1\times N_v \times r_2$.
$G_\rom{3}$ and $H_\rom{3}$ are matrices of sizes $r_1\times r_1$ and $r_2\times r_2$, respectively, while 
$T^{(2)}$ is a tridiagonal matrix of size $N_v \times N_v$. Again, the discussion of how to solve this equation is deferred to the next section.

Finally, in \textbf{substep 5}, the unknown tensor train $\tilde{\bf}^j$ is in the form $(\rom{5})$, denoted by $\llbracket P^{(1)}, P^{(2)}, \tilde{C}^{(3)}\rrbracket$, and only the third core $\tilde{C}^{(3)}$ is updated. This update is performed through the projection of \eqref{eq_forward_scheme_transformed} onto the cores $P^{(1)}$ and $P^{(2)}$. Diagrammatically, this means that the whole equation is contracted with $P^{(1)}$ and $P^{(2)}$:
\begin{equation}
\label{diag_substep5_projection}
\underbrace{
\begin{tikzpicture}[baseline=-26]
\draw[rounded corners=2pt, dashed, gray] (-0.5,1) rectangle (1.7,-3);
\draw (0,0) -- (1.2,0);
\draw (1.2,0) -- (2.4,0);
\draw (0,0) -- (0,-2);
\draw (1.2,0) -- (1.2,-2);
\draw (2.4,0) -- (2.4,-1.6);
\draw[fill=white,line width = 1.2] (2.4,0) circle (0.2);
\draw[fill=white,line width=1.2] (1.4,0) -- (1.1,-0.2) -- (1.1,0.2) -- (1.4,0);
\draw[fill=white,line width=1.2] (0.2,0) -- (-0.1,-0.2) -- (-0.1,0.2) -- (0.2,0);
\node at (2.4,0.5) {$\tilde{C}^{(3)}$};
\node at (1.2,0.5) {${P}^{(2)}$};
\node at (0,0.5) {${P}^{(1)}$};
\draw [fill=black] (0,-0.8) -- (-0.2,-1) -- (0,-1.2) -- (0.2,-1);
\node at (0.6,-1) {$T^{(1)}$};
\draw (2,-2) -- (0,-2);
\draw[fill=white,line width=1.2] (1.4,-2) -- (1.1,-2.2) -- (1.1,-1.8) -- (1.4,-2);
\draw[fill=white,line width=1.2] (0.2,-2) -- (-0.1,-2.2) -- (-0.1,-1.8) -- (0.2,-2);
\node at (1.2,-2.5) {${P}^{(2)}$};
\node at (0,-2.5) {${P}^{(1)}$};
\end{tikzpicture}
+
\begin{tikzpicture}[baseline=-26]
\draw[rounded corners=2pt, dashed, gray] (-0.5,1) rectangle (1.7,-3);
\draw (0,0) -- (1.2,0);
\draw (1.2,0) -- (2.4,0);
\draw (0,0) -- (0,-2);
\draw (1.2,0) -- (1.2,-2);
\draw (2.4,0) -- (2.4,-1.6);
\draw[fill=white,line width = 1.2] (2.4,0) circle (0.2);
\draw[fill=white,line width=1.2] (1.4,0) -- (1.1,-0.2) -- (1.1,0.2) -- (1.4,0);
\draw[fill=white,line width=1.2] (0.2,0) -- (-0.1,-0.2) -- (-0.1,0.2) -- (0.2,0);
\node at (2.4,0.5) {$\tilde{C}^{(3)}$};
\node at (1.2,0.5) {${P}^{(2)}$};
\node at (0,0.5) {${P}^{(1)}$};
\draw [fill=black] (1.2,-0.8) -- (1,-1) -- (1.2,-1.2) -- (1.4,-1);
\node at (1.8,-1) {$T^{(2)}$};
\draw (2,-2) -- (0,-2);
\draw[fill=white,line width=1.2] (1.4,-2) -- (1.1,-2.2) -- (1.1,-1.8) -- (1.4,-2);
\draw[fill=white,line width=1.2] (0.2,-2) -- (-0.1,-2.2) -- (-0.1,-1.8) -- (0.2,-2);
\node at (1.2,-2.5) {${P}^{(2)}$};
\node at (0,-2.5) {${P}^{(1)}$};
\end{tikzpicture}
}_{\tilde{C}^{(3)}\times_1 G_\rom{5}}
+
\underbrace{
\begin{tikzpicture}[baseline=-26]
\draw[rounded corners=2pt, dashed, gray] (-0.5,1) rectangle (1.7,-3);
\draw (0,0) -- (1.2,0);
\draw (1.2,0) -- (2.4,0);
\draw (0,0) -- (0,-2);
\draw (1.2,0) -- (1.2,-2);
\draw (2.4,0) -- (2.4,-1.6);
\draw[fill=white,line width = 1.2] (2.4,0) circle (0.2);
\draw[fill=white,line width=1.2] (1.4,0) -- (1.1,-0.2) -- (1.1,0.2) -- (1.4,0);
\draw[fill=white,line width=1.2] (0.2,0) -- (-0.1,-0.2) -- (-0.1,0.2) -- (0.2,0);
\node at (2.4,0.5) {$\tilde{C}^{(3)}$};
\node at (1.2,0.5) {${P}^{(2)}$};
\node at (0,0.5) {${P}^{(1)}$};
\draw [fill=black] (2.4,-0.8) -- (2.2,-1) -- (2.4,-1.2) -- (2.6,-1);
\node at (1.8,-1) {$T^{(3)}$};
\draw (2,-2) -- (0,-2);
\draw[fill=white,line width=1.2] (1.4,-2) -- (1.1,-2.2) -- (1.1,-1.8) -- (1.4,-2);
\draw[fill=white,line width=1.2] (0.2,-2) -- (-0.1,-2.2) -- (-0.1,-1.8) -- (0.2,-2);
\node at (1.2,-2.5) {${P}^{(2)}$};
\node at (0,-2.5) {${P}^{(1)}$};
\end{tikzpicture}
}_{\tilde{C}^{(3)} \times_2 T^{(3)}}
=
\underbrace{
\begin{tikzpicture}[baseline=-15]
\draw (0,0) -- (1.2,0);
\draw (1.2,0) -- (2.4,0);
\draw (1.2,-1.2) -- (2,-1.2);
\draw (0,-1.2) -- (1.2,-1.2);
\draw (0,0) -- (0,-1.2);
\draw (1.2,0) -- (1.2,-1.2);
\draw (2.4,0) -- (2.4,-0.8);
\draw[fill=white,line width = 1.2] (0,0) circle (0.2);
\draw[fill=white,line width = 1.2] (1.2,0) circle (0.2);
\draw[fill=white,line width = 1.2] (2.4,0) circle (0.2);
\draw[fill=white,line width=1.2] (0.2,-1.2) -- (-0.1,-1.4) -- (-0.1,-1) -- (0.2,-1.2);
\draw[fill=white,line width=1.2] (1.4,-1.2) -- (1.1,-1.4) -- (1.1,-1) -- (1.4,-1.2);
\node at (0,0.5) {${K}^{(1)}$};
\node at (1.2,0.5) {${K}^{(2)}$};
\node at (2.4,0.5) {${K}^{(3)}$};
\node at (0,-1.7) {${P}^{(1)}$};
\node at (1.2,-1.7) {${P}^{(2)}$};
\end{tikzpicture}
}_{R_\rom{5}}
\end{equation}
Similar to substep 1, this yields a matrix Sylvester equation:
\begin{equation}
\label{eq_Sylvester_substep5}
    \tilde{C}^{(3)}\times_1 G_\rom{5} + 
    \tilde{C}^{(3)}\times_2 T^{(3)} = R_\rom{5},
\end{equation}
where the unknown $\tilde{C}^{(3)}$ and the right-hand side $R_\rom{5}$ are matrices of size $r_2\times N_v$;
$T^{(3)}$ is a tridiagonal matrix of size $N_v\times N_v$, while $G_\rom{5}$ is a matrix of size $r_2\times r_2$.

\begin{remark}
Although the discussion in this paper focuses on the Fokker--Planck operator, the implicit method presented here can be naturally extended to other diffusion-type equations. In particular, if $\mathcal{M}[f]$ in \eqref{eq_Fokker_Planck_collision} is fixed to be the constant $1$, the operator reduces to the Laplacian. In this case, the method can be applied directly to the heat equation. The method can also be applied to the variable-coefficient diffusion equation of the form:
\begin{equation}
    \partial_t f(t,\bv) = \nabla_\bv \cdot \left(A(\bv)
    \nabla_\bv f
    \right),
    \text{~with~}
    A(\bv) = \operatorname{diag}(a_1(v_1),a_2(v_2),a_3(v_3)).
\end{equation}
For more general variable-coefficient diffusion equations, however, the method requires careful adaptation depending on the structure of the coefficients. We defer these extensions to future work.

\end{remark}

\section{Sylvester solvers}
\label{sec_Sylvester_Solvers}
The scheme proposed in the previous section involves solving the Sylvester equations \eqref{eq_Sylvester_substep1}, \eqref{eq_Sylvester_substep3}, and \eqref{eq_Sylvester_substep5}.
The Sylvester equation is an important problem in numerical linear algebra and has been extensively studied; see, for example, the review \cite{simoncini2016computational}. 
In this section, we show that a direct application of off-the-shelf Sylvester solvers to our problems results in an undesirably high computational complexity of $\mathcal{O}(N_v^3)$. By exploiting the special structure of the matrices, we propose a tailored direct Sylvester solver that achieves linear complexity with respect to $N_v$.
%

%
We begin with the matrix Sylvester equation appearing in \textbf{substeps 1} and \textbf{5}:
\begin{equation}
\label{eq_matrix_sylvester}
    X \times_1 L_1 + X \times_2 L_2 = R.
\end{equation}
By the definition of the operator $\times_n$, Eq.~\eqref{eq_matrix_sylvester} is equivalent to the standard matrix form 
\begin{equation}
\label{eq_matrix_sylvester1}
L_1^\top X + X L_2 = R.
\end{equation}
In \textbf{substep 1} (Eq.~\eqref{eq_Sylvester_substep1}), $L_1$ is a matrix of size $N_v\times N_v$ and $L_2$ is a matrix of size $r_1\times r_1$. In \textbf{substep 5} (Eq.~\eqref{eq_Sylvester_substep5}), $L_1$ is a matrix of size $r_2\times r_2$ and $L_2$ is a matrix of size $N_v\times N_v$. A classical approach for solving these equations is the Bartels--Stewart algorithm \cite{BS72}, which 
computes the Schur decompositions of both matrices $L_1$ and $L_2$, and then transforms the equation into a lower/upper triangular form that can be solved explicitly element by element.
This algorithm is implemented in standard dense linear algebra software packages, such as the Sylvester solvers in SciPy and MATLAB.
The computational cost is dominated by the Schur decompositions. In our setting, since one of the matrices has the large size $N_v\times N_v$ (typically $N_v\gg r_1,r_2$), the overall computational complexity becomes 
$\mathcal{O}(N_v^3)$, which is prohibitively expensive. In fact, this would make the entire low-rank algorithm as expensive as the full-tensor method.

On the other hand, we observe that our problems exhibit a special structure that can be further exploited to design a more efficient solver. Let us take \eqref{eq_Sylvester_substep1} as an example, where $L_1=T^{(1)}\in \mathbb{R}^{N_v\times N_v}$ and $L_2=H_\rom{1}\in \mathbb{R}^{r_1\times r_1}$ (the case of \eqref{eq_Sylvester_substep5} is analogous). Although $L_1$ is large, it is tridiagonal. Applying a Schur decompositions to $L_1$ is therefore inefficient and would also destroy its sparsity pattern. Motivated by this observation, we compute the Schur decomposition only for the small matrix $L_2 = U_2 W_2 U_2^*$, where $U_2$ is a unitary matrix, $U_2^*$ is its conjugate transpose, and $W_2$ is an upper triangular matrix. 
Define 
$Z=XU_2$, $E=RU_2$, the equation \eqref{eq_matrix_sylvester1} then becomes
\begin{equation}
\label{eq_matrix_Sylvester_reduced}
    L_1^\top Z + Z W_2 = E.
\end{equation}
The columns of $Z$ can then be solved sequentially using the tridiagonal Thomas algorithm.
More specifically, the $j$th column of $Z$, denoted by $Z_{:,j}$, satisfies 
\begin{equation}
\label{eq_tridiagonal}
    (L_1^\top + (W_2)_{jj} I) Z_{:,j} = E_{:,j} - \sum_{k=1}^{j-1}(W_2)_{kj} Z_{:,k},
\end{equation}
where the coefficient matrix $L_1^\top + (W_2)_{jj} I$ remains tridiagonal. Consequently, the matrix $Z$ can be computed column by column from left to right. Finally, the matrix $X$ is recovered via $X = ZU_2^*$.

We estimate the computational cost of the above proposed algorithm. 
The Schur decomposition of $L_2$ requires $\mathcal{O}(r_1^3)$ operations. 
Forming the intermediate matrix $E$ and recovering the solution $X$ from $Z$ both incur a cost of $\mathcal{O}(N_v r_1^2)$. Solving \eqref{eq_tridiagonal} using the Thomas algorithm for each $j$ costs $\mathcal{O}(N_v)$, while forming the right-hand side costs $O(N_v j)$. Looping over $r_1$ columns of $Z$, this part costs $\mathcal{O}(N_v r_1^2)$. Altogether, since $N_v\gg r_1$, the overall computational complexity is dominated by $\mathcal{O}(N_v r_1^2)$, which scales linearly with respect to $N_v$ and is therefore well suited for large-scale problems.

We now turn to the tensor Sylvester equation appearing in \textbf{substep 3} (Eq.~\eqref{eq_Sylvester_substep3}):
\begin{equation}
\label{tensor-sylvester}
    X \times_1 L_1 + X \times_2 L_2 + X \times_3 L_3 = R,
\end{equation}
where $L_1=G_{\rom{3}}\in \mathbb{R}^{r_1\times r_1}$ and $L_3=H_{\rom{3}}\in \mathbb{R}^{r_2\times r_2}$ are two small dense  matrices, while $L_2=T^{(2)}\in \mathbb{R}^{N_v\times N_v}$ is a large tridiagonal matrix. One may use the Bartels--Stewart algorithm to solve this tensor equation, which again leads to an $\mathcal{O}(N_v^3)$  complexity.

Following a similar strategy as above, we compute the Schur decompositions only for $L_1$ and $L_3$: $L_1=U_1 W_1 U_1^*$, $L_3=U_3 W_3 U_3^*$. Define $Z = (X \times_1 U_1) \times_3 U_3$ and $E = (R \times_1 U_1) \times_3 U_3$,  \eqref{tensor-sylvester} is then reduced to
\begin{equation}
\label{eq_tensor_Sylvester_reduced}
    Z \times_1 W_1 + Z \times_2 L_2 + Z \times_3 W_3 = E.
\end{equation}
This equation can be solved by considering each fiber of $Z$ along the second mode. 
Fix indices $(i,j)$ and denote by $Z_{i,:,j} \in \mathbb{C}^{N_v}$ the corresponding vector, and similarly $E_{i,:,j}$. 
The vector $Z_{i,:,j}$ satisfies
\begin{equation}
\label{eq_tensor_Sylvester_slice}
    (L_2 + (W_1)_{ii} I
    + (W_3)_{jj} I )Z_{i,:,j}
    = E_{i,:,j} - \sum_{p=1}^{i-1} (W_1)_{pi} Z_{p,:,j}
    - \sum_{q=1}^{j-1} (W_3)_{qj} Z_{i,:,q}.
\end{equation}
Here, the right-hand side depends only on entries that have already been computed, which makes the recursion well-defined.
The resulting systems retain a tridiagonal structure in the $L_2$ direction and can therefore be solved efficiently using the Thomas algorithm. Once $Z$ has been computed, the original tensor $X$ is recovered by 
$X = (Z \times_1 U_1^*) \times_3 U_3^*$.

For this algorithm, the two Schur decompositions require $\mathcal{O}(r_1^3 + r_2^3)$ operations. 
Forming the intermediate tensor $E$ and transforming $Z$ back to $X$ involve contractions along the first and third modes and cost $\mathcal{O}(r_1^2 N_v r_2 + r_1 N_v r_2^2)$. 
For each fiber equation 
\eqref{eq_tensor_Sylvester_slice}, forming the right-hand side requires $\mathcal{O}(N_v(r_1+r_2))$ operations, while the tridiagonal solve itself costs $\mathcal{O}(N_v)$. Since there are $r_1r_2$ such fiber equations, the total cost of this part is $\mathcal{O}(r_1 r_2 N_v (r_1+r_2))$. Combining everything, the overall computational complexity is $\mathcal{O}(r_1^2 N_v r_2 + r_1 N_v r_2^2)$, which is again linear with respect to the large dimension $N_v$.

In the above discussion, complex Schur decompositions are employed, which introduce complex-valued intermediate quantities. 
In the absence of rounding errors, however, the resulting solution $X$ is guaranteed to be real-valued, either as a matrix or as a tensor. 
In practice, due to numerical errors, we take the real part of $X$ as the solution of the Sylvester equation.
If one wishes to avoid complex arithmetic, one may instead employ the real Schur decomposition. 
In this case, the resulting matrix is quasi-triangular and may contain $2 \times 2$ blocks on the diagonal. 
Consequently, the corresponding unknowns must be solved in coupled form: two columns of $Z$ in \eqref{eq_matrix_Sylvester_reduced}, and two or four fibers of $Z$ in \eqref{eq_tensor_Sylvester_reduced}, depending on the block structure.

\section{Memory requirement and computational complexity}
\label{sec_cost}

In this section, we estimate the computational complexity of the proposed low-rank algorithm in Section~\ref{sec_integrator}, together with the tailored Sylvester solvers presented in Section~\ref{sec_Sylvester_Solvers}.

First, in terms of memory, the full tensor method requires $\mathcal{O}(N_xN_v^3)$ storage, whereas the low-rank algorithm requires only $\mathcal{O}(r^2N_xN_v)$, where $r=\max\{r_1,r_2\}$ denotes the maximal TT-rank.

In terms of computational complexity, the full tensor method requires at least $\mathcal{O}(N_x N_v^3)$ operations per time step, and often more, typically $\mathcal{O}(KN_x N_v^3)$, where $K$ is the number of iterations, since the implicit Fokker--Planck operator generally requires iterative solves of a large three-dimensional tensor Sylvester equation at each spatial grid point.

For the implicit dynamical TT algorithm proposed in this paper, at each time step we first solve \eqref{eq_moments} to obtain the macroscopic quantities and hence the Maxwellian. Since the solution is represented in TT format, this step can be carried out efficiently with a cost of $\mathcal{O}(r^2 N_x N_v)$ \cite{oseledets2011tensor,wang2026dynamical,wang2026accelerated}.
The subsequent five-substep procedure consists of repeated contractions between tensor cores and matrices. 
Each substep involves operations with complexity $\mathcal{O}(r^3 N_v)$ as discussed in \cite{wang2026dynamical}. 
In addition, solving the matrix Sylvester equations in substeps 1 and 5 requires $\mathcal{O}(r^2N_v)$ operations, while solving the tensor Sylvester equations in substep 3 requires $\mathcal{O}(r^3N_v)$ operations.
Since these computations are performed independently at each spatial grid point, the total cost of the five-substep process scales as $\mathcal{O}(r^3 N_x N_v)$.

Combining the above contributions, the overall computational complexity of the proposed algorithm is $\mathcal{O}(r^3 N_x N_v)$ per time step. We note that this cost is close to optimal, as it is of the same order as that of the explicit dynamical TT method in \cite{wang2026dynamical}.

\section{Numerical examples}
\label{sec_numerical_results}
In this section, we apply the algorithm introduced in the previous sections to several kinetic equations with Fokker--Planck collision term and demonstrate its performance. All the experiments were conducted on a MacBook Pro equipped with an Apple M4 CPU.

\subsection{Spatially homogeneous Fokker--Planck equation}
We first consider the spatially homogeneous (linear) Fokker--Planck equation, which reads
\begin{equation}
\label{homo}
    \partial_t f(t,\bv) = \nabla_\bv \cdot\left( M(\bv) \nabla_\bv \left(\frac{f}{M(\bv)}\right)
    \right),
\end{equation}
where the Maxwellian $M(\bv)$ is given by
\begin{equation}
    M(\bv) = \frac{1}{(2\pi)^{3/2}}\exp\left(-\frac{\vert \bv \vert^2}{2}\right).
\end{equation}
If the initial condition is chosen as
\begin{equation}
    f_0(\bv) = \frac{1}{(2\pi(1-\e^{-1}))^{3/2}} \exp\left(-\frac{\vert \bv\vert^2}{2(1-\e^{-1})}\right),
\end{equation}
then (\ref{homo}) admits an exact solution
\begin{equation}
    f_{\text{exact}}(t,\bv) = \frac{1}{(2\pi(1-\e^{-(1+2 t)}))^{3/2}}
    \exp\left(-\frac{\vert\bv\vert^2}{2(1-\e^{-(1+2 t)})}\right).
\end{equation}
\begin{figure}
    \centering
    \includegraphics[width=0.55\linewidth]{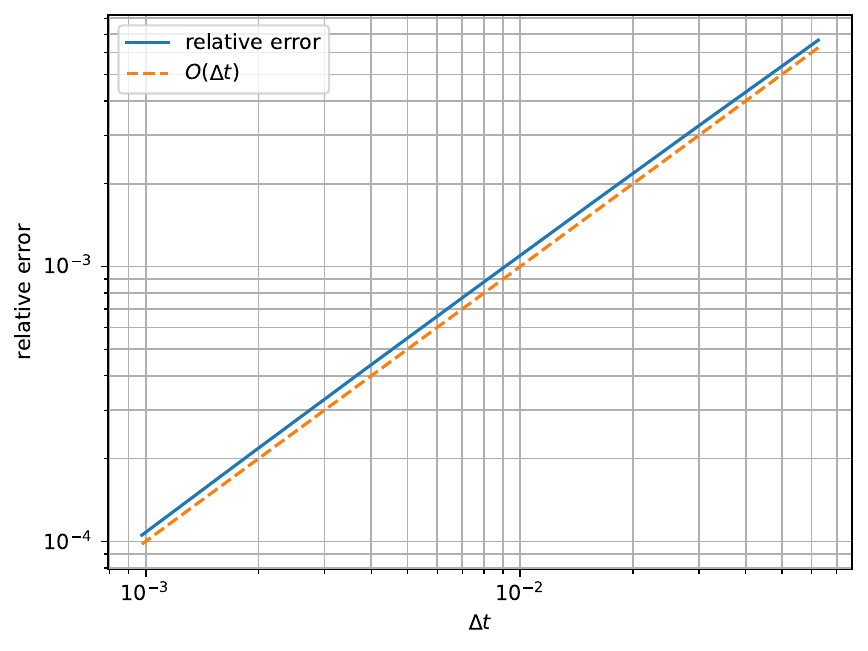}
    \caption{Spatially homogeneous Fokker--Planck equation: relative error with different $\Delta t$.}
    \label{fig_linear_Fokker_Planck}
\end{figure}

We solve (\ref{homo}) up to time $t=1$ using a fixed TT-rank of $(5,5)$. 
The velocity domain is taken as $[v_{\min}, v_{\max}]^3 = [-8,8]^3$ with a fixed mesh size $\Delta v = \frac{1}{32}$ and grid points $
    v_k = v_{\min} + \left(k-\frac{1}{2}\right)\dv$, $k=1,\cdots,N_v$.
The time step $\Delta t$ varies from $\frac{1}{16}$ to $\frac{1}{1024}$. We note that these time steps are chosen independently of $\Delta v$, precisely due to the advantage of the implicit scheme.

The relative error, defined by $\dfrac{\Vert \bf - \bf_{\mathrm{exact}} \Vert_F}{\Vert \bf_{\mathrm{exact}} \Vert_F}$, is shown in \Cref{fig_linear_Fokker_Planck}.
The results clearly indicate that the numerical error scales linearly with respect to the time step $\Delta t$.

\begin{figure}
    \centering
    \includegraphics[width=0.55\linewidth]{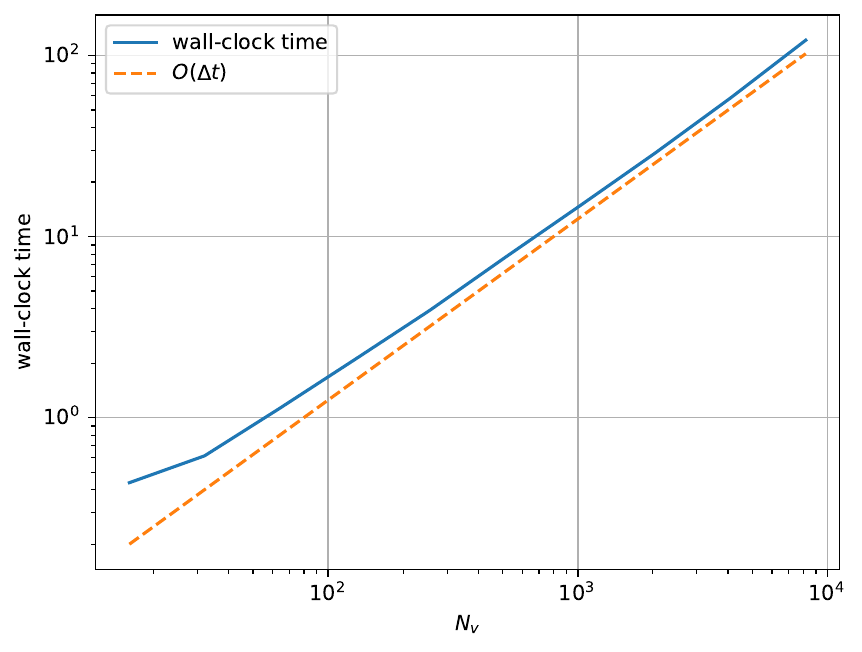}
    \caption{Spatially homogeneous Fokker--Planck equation: wall-clock time for different $N_v$.}
    \label{fig_time}
\end{figure}

To demonstrate the linear scaling of the computational complexity with respect to $N_v$ as analyzed in \Cref{sec_cost}, 
we perform numerical experiments with varying $\Delta v$ while fixing $\dt = \frac{1}{512}$. 
The corresponding wall-clock times are reported in \Cref{fig_time}, which 
clearly exhibits linear scaling with respect to $N_v$.

\subsection{Spatially inhomogeneous Fokker--Planck equation}
We then consider the 1D3V spatially inhomogeneous Fokker--Planck equation:
\begin{equation}
    \partial_t f(t,x,\bv) + v^{(1)}\partial_{x} f(t,x,\bv) = \eta \, T \,
    \nabla_\bv\cdot \left(
        \mathcal{M}[f] \nabla_\bv \left(\frac{f}{\mathcal{M}[f]}\right)
    \right),
\end{equation}
with spatial domain $x\in[0,1]$ and periodic boundary conditions.
The initial condition is given by 
\begin{equation}
    f_0(x,\bv) = \frac{n_0(x)}{(2\pi T_0(x))^{3/2}} \exp\left(-\frac{\vert \bv - \bu_0 \vert^2}{2T_0(x)}\right),
\end{equation}
with
\begin{equation}
    n_0(x) = \frac{2+\sin(2\pi x)}{3}, \quad
    \bu_0=(0.2,0,0)^\top, \quad
    T_0(x) = \frac{3+\cos(2\pi x)}{4}.
\end{equation}
In the numerical experiment, we choose $N_x=64$, $\dx = \dfrac{1}{N_x}$, and uniform grid points 
   $ x_j = \left(j-\frac{1}{2}\right) \dx$, $j = 1,\cdots,N_x$.
The velocity domain is truncated to $[v_{\min},v_{\max}]^3=[-6,6]^3$ with $N_v = 64$ uniform grid points in each direction.
The transport term $v^{(1)} \partial_x f(t,x,\bv)$ is discretized by the second-order upwind scheme
\begin{equation}
\label{eq_transport_term_upwind}
    \left(
        v^{(1)} \bD_x^\up \bf
    \right)^{j,n}
    = \frac{v_{k_1}^+}{2\dx}
    (3\bf^{j,n}-4\bf^{j-1,n}+\bf^{j-2,n})
    - \frac{v_{k_1}^-}{2\dx}
    (-\bf^{j+2,n}+4\bf^{j+1,n}-3\bf^{j,n}),
\end{equation}
with $v_{k_1}^+ = \max\{v_{k_1},0\}$,
$v_{k_1}^- = \max\{-v_{k_1},0\}$.
When $\bf^{j,n}$ is in the TT form, the transport term $\left(v^{(1)} \bD_x^\up \bf\right)^{j,n}$ can also be represented in TT form.

\begin{figure}[t]
     \centering
     \begin{subfigure}[b]{0.355\textwidth}
         \centering
         \includegraphics[width=\textwidth]{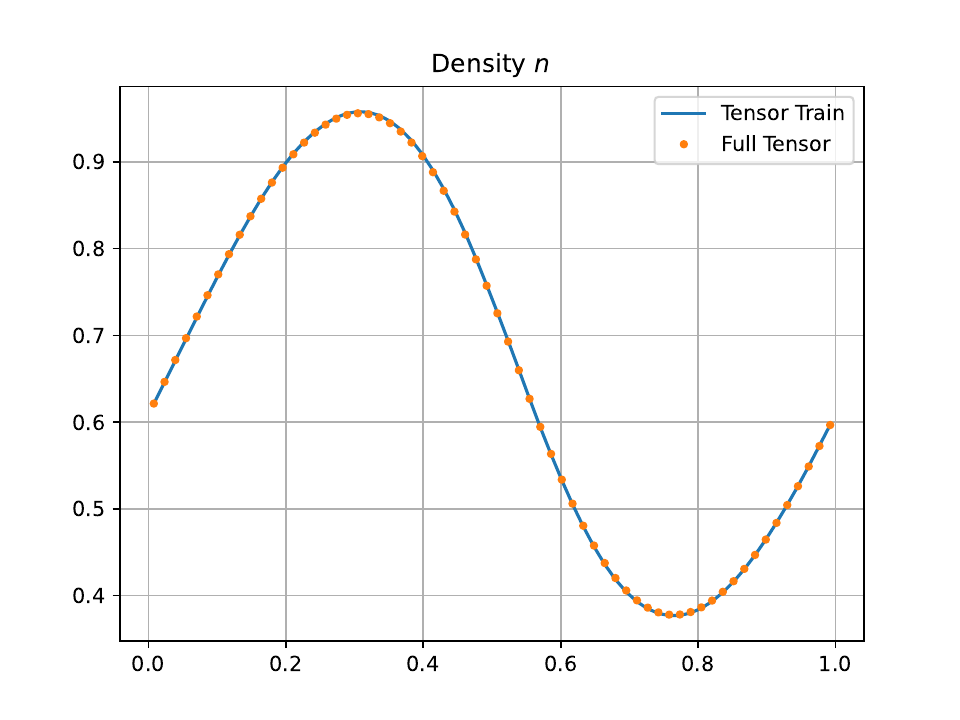}
         \caption{Density at $t=0.1$.}
     \end{subfigure}
     \hspace{-22pt}
     \begin{subfigure}[b]{0.355\textwidth}
         \centering
         \includegraphics[width=\textwidth]{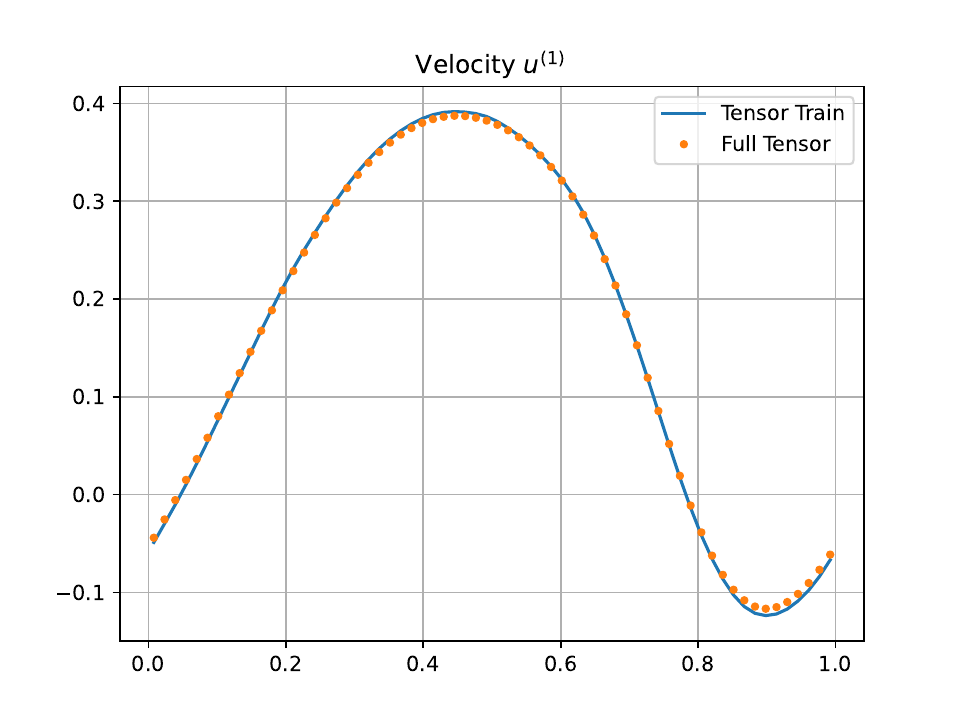}
         \caption{Velocity at $t=0.1$.}
     \end{subfigure}
     \hspace{-22pt}
     \begin{subfigure}[b]{0.355\textwidth}
         \centering
         \includegraphics[width=\textwidth]{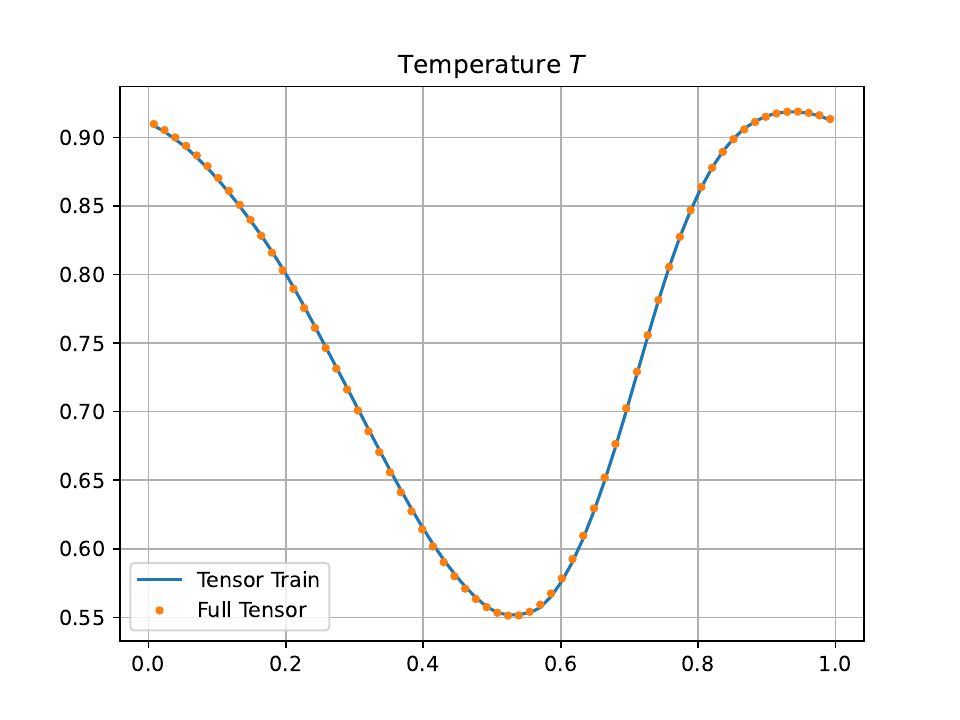}
         \caption{Temperature at $t=0.1$.}
     \end{subfigure}
        \caption{Spatially inhomogeneous Fokker--Planck equation: density, bulk velocity, and temperature for $\eta = 1$.}
        \label{fig_inhomo_FP_eta1}
\end{figure}

\begin{figure}[t]
     \centering
     \begin{subfigure}[b]{0.355\textwidth}
         \centering
         \includegraphics[width=\textwidth]{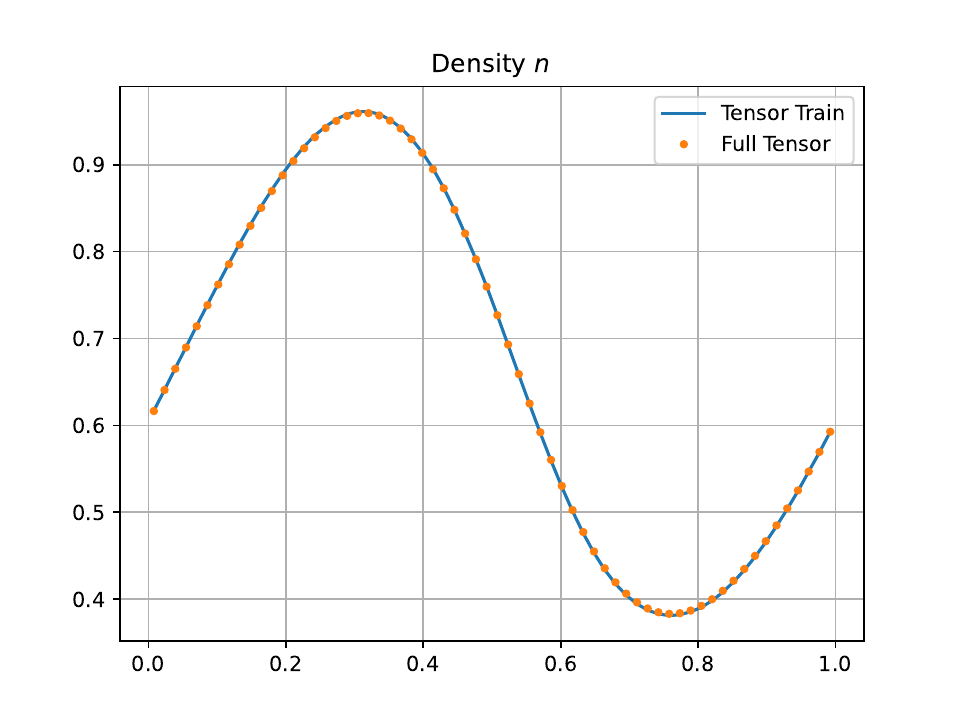}
         \caption{Density at $t=0.1$.}
     \end{subfigure}
     \hspace{-22pt}
     \begin{subfigure}[b]{0.355\textwidth}
         \centering
         \includegraphics[width=\textwidth]{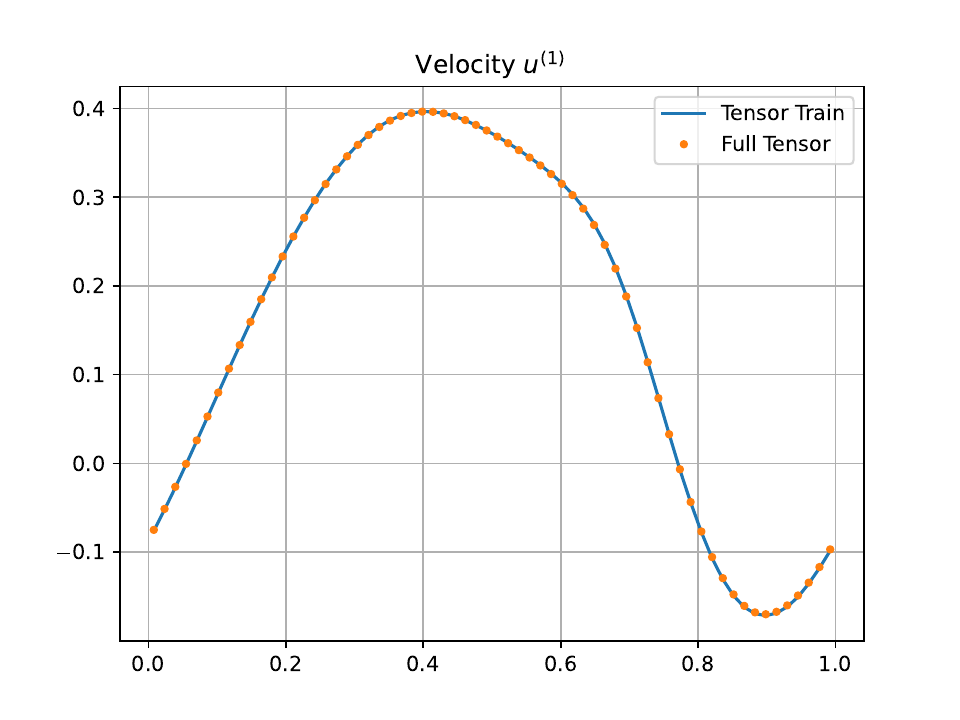}
         \caption{Velocity at $t=0.1$.}
     \end{subfigure}
     \hspace{-22pt}
     \begin{subfigure}[b]{0.355\textwidth}
         \centering
         \includegraphics[width=\textwidth]{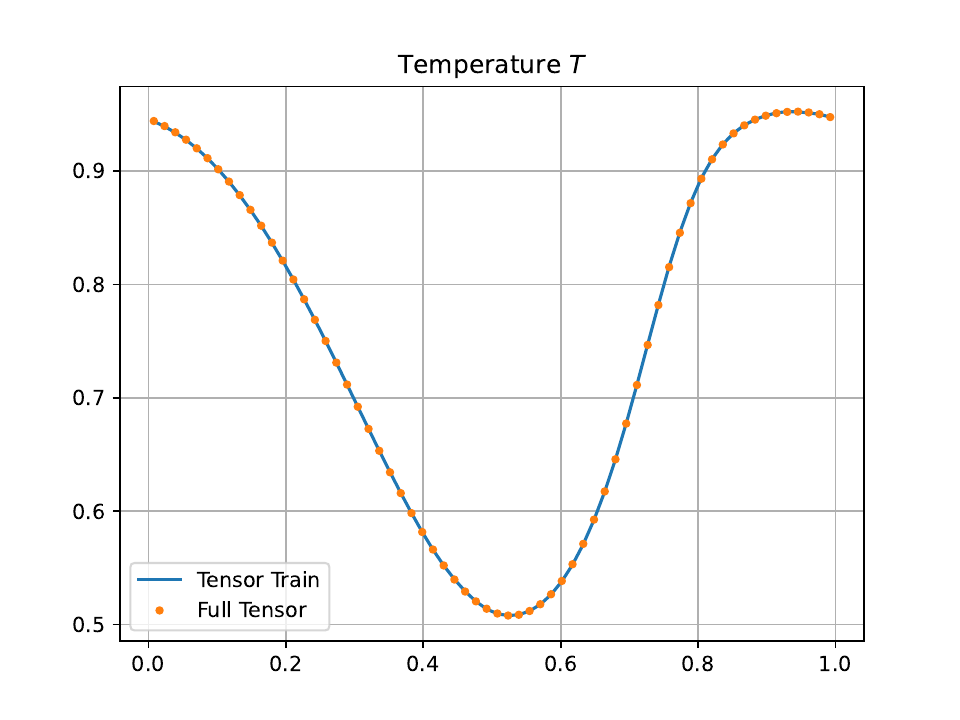}
         \caption{Temperature at $t=0.1$.}
     \end{subfigure}
        \caption{Spatially inhomogeneous Fokker--Planck equation: density, bulk velocity, and temperature for $\eta = 10^6$.}
        \label{fig_inhomo_FP_eta1000000}
\end{figure}

We set $\eta = 1$ and $\eta = 10^6$, corresponding to the transitional and stiff regimes, respectively.
The equation is solved up to $t=0.1$ with a time step $\dt = 0.001$. The TT-rank is fixed as $(5,5)$. The macroscopic quantities $n$, $\bu$, and $T$, defined in (\ref{eq_macroscopic_quantities}), are presented, along with
the results from the full tensor method using the first-order IMEX scheme.
From \Cref{fig_inhomo_FP_eta1,fig_inhomo_FP_eta1000000}, we observe that the results produced by our method agree well with those obtained from the full tensor scheme. 
At the same time, our method is significantly more efficient in terms of wall-clock time. 
For instance, in this example with 100 time steps, our method completes in less than one minute, whereas the full tensor scheme requires more than one hour.

Regarding the admissible time step, in a fully explicit scheme, the time step must satisfy the restriction
$\dt \lesssim \min\left(\frac{\dx}{\max |v^{(1)}|}, \frac{\dv^2}{6\eta}\right)$ to ensure stability,
where the first constraint arises from the transport term, while the second is imposed by the collision term \cite{wang2026dynamical}. 
In the proposed method, the second restriction is removed due to the implicit treatment of the collision term. 
As a result, significantly larger time steps can be employed, particularly in the stiff regime. 
For example, in our experiment, when $\eta = 10^6$, a fully explicit scheme such as that in \cite{wang2026dynamical} requires $\dt \lesssim 4 \times 10^{-8}$, whereas the present method allows $\dt \lesssim 2.6 \times 10^{-3}$. 
This corresponds to a difference of several orders of magnitude, demonstrating that the implicit treatment leads to a substantial acceleration of the numerical simulation.


\subsection{Vlasov--Ampère--Fokker--Planck (VAFP) equation}
In this section, we apply our method to the Vlasov--Fokker--Planck equation in the 1D3V setting.
The equation reads
\begin{equation}
    \partial_t f+v^{(1)} \partial_x f-E^{(1)} \partial_{v^{(1)}} f=\eta T \nabla_{\mathbf{v}} \cdot\left(\mathcal{M}[f] \nabla_{\mathbf{v}}\left(\frac{f}{\mathcal{M}[f]}\right)\right),
\end{equation}
where the electric field follows the Ampère's law
\begin{equation}
    \partial_t E^{(1)}(t,x) = -J^{(1)}(t,x), \text{~with~}
    J^{(1)}(t,x) 
    = -\int_{\mathbb{R}^3} v^{(1)} f(t,x,\bv) \rd \bv.
\end{equation}
The discretization of the transport term is the same as (\ref{eq_transport_term_upwind}). The force term is also discretizd by the second-order upwind:
\begin{equation}
    (E^{(1)} \mathbf{D}_{v^{(1)}}^{\mathrm{up}} \mathbf{f})_{k_1 k_2 k_3}^j=\frac{(E_j^{(1)})^{+}(-f_{k_1+2, k_2 k_3}^j+4 f_{k_1+1, k_2 k_3}^j-3 f_{k_1 k_2 k_3}^j)-(E_j^{(1)})^{-}(3 f_{k_1 k_2 k_3}^j-4 f_{k_1-1, k_2 k_3}^j+f_{k_1-2, k_2 k_3}^j)}{2 \Delta v},
\end{equation}
with $(E_j^{(1)})^{+}=\max\{(E_j^{(1)}),0\}$ and $(E_j^{(1)})^{-}=\max\{-(E_j^{(1)}),0\}$. We apply the zero boundary condition in the $v^{(1)}$ direction by setting $f_{0,k_2k_3}^j=f_{-1,k_2k_3}^j=f_{N_v+1,k_2k_3}^j=f_{N_v+2,k_2k_3}^j=0$.

We consider the linear Landau damping with  initial condition
\begin{equation}
\label{eq_damping_initial}
    f_0(x,\bv) = \frac{1}{(2\pi)^{3/2}} \bigl(1 + A \cos(\kappa x)\bigr) \e^{-\vert \bv \vert^2/2}.
\end{equation}
The initial electric field is determined from Gauss's law:
\begin{equation}
    \frac{\partial E^{(1)}}{\partial x} = \rho - \rho_i,
\end{equation}
where $\displaystyle \rho = -\int_{\mathbb{R}^3} f \, \dd \bv$ denotes the charge density, and $\rho_i$ is a uniform background density introduced to neutralize the system. 
For the initial condition (\ref{eq_damping_initial}), the corresponding electric field is given by
\begin{equation}
    E^{(1)}(0,x) = -\frac{A}{\kappa} \sin(\kappa x).
\end{equation}

The physical parameters are chosen as $A = 0.001$ and $\kappa = 0.5$. 
The spatial domain is $x \in [0, 2\pi/\kappa] = [0, 4\pi]$ with $N_x = 128$ grid points and periodic boundary condition, while the velocity domain is truncated to $[v_{\min},v_{\max}]^3=[-9,9]^3$ with $N_v = 128$ points in each direction.
Again since the Fokker--Planck operator is treated implicitly, the time step is no longer constrained by the collision term, in contrast to the fully explicit scheme in \cite{wang2026dynamical}.
The time step is therefore selected as
\begin{equation}
\label{eq_vlasov_time_step}
    \dt = 0.1 \min\left\{\frac{\dx}{\max \left\vert v^{(1)}\right\vert}, \frac{\dv}{\max \left\vert E^{(1)}(0,x) \right\vert}\right\}.
\end{equation}
The TT-rank is fixed as $(5,5)$ throughout the simulation.
To investigate the effect of collisions, we consider $\eta = 0$, $1$, and $2$.
With the chosen parameter settings, the time step is fixed as $\dt \approx 0.00109$. 
For comparison, the time steps required by the fully explicit scheme in \cite{wang2026dynamical} are $0.00109$ for $\eta = 0$, $0.00033$ for $\eta = 1$, and $0.00016$ for $\eta = 2$. 
This clearly illustrates that, as the collision strength increases, the explicit scheme becomes increasingly restricted, whereas the present method allows a relatively large time step.
The electric energy is defined as%
\begin{equation}
    \mathcal{E}(t) = \frac{1}{2} \int_0^{2\pi/\kappa} \bigl(E^{(1)}(t,x)\bigr)^2 \dd x
    \approx \frac{1}{2} \sum_{j=1}^{N_x} \left(E^{(1)}_j\right)^2 \dx.
\end{equation}
\begin{figure}
    \centering
    \includegraphics[width=0.55\linewidth]{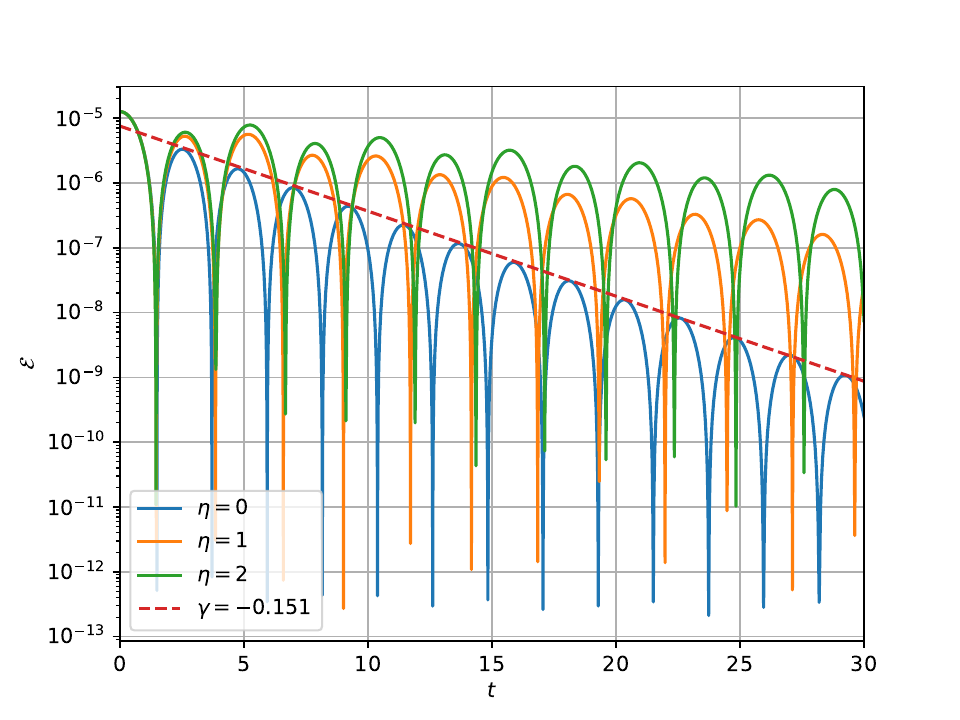}
    \caption{Landau damping. Evolution of electric energy for different collision strengths.}
    \label{fig_Vlasov_LD}
\end{figure}%
The evolution of the electric energy is shown in Figure \ref{fig_Vlasov_LD}.
In the collisionless case ($\eta = 0$), we observe a clear damping rate of $\gamma = -0.151$. 
This result is consistent with the fully explicit scheme in \cite{wang2026dynamical} and agrees well with the theoretical prediction ($\gamma = -0.153$) \cite{cheng1976integration}. 
In the collisional case ($\eta = 1,2$), however, a single dominant damping rate is no longer observed. 
Instead, the solution exhibits an oscillatory damping behavior characterized by alternating high and low peaks. 
This pattern indicates that the dynamics are governed by the superposition of multiple damped modes with comparable amplitudes but distinct frequencies, resulting in a beating-like modulation in time. 
Such behavior has also been reported in previous studies \cite{ng2004complete,black2013discrete,jorge2019linear,banik2024relaxation,ye2024energy,coughlin2024robust}.
For moderate collisional strengths, both kinetic effects and collisional relaxation play important roles, leading to the coexistence of multiple discrete modes.

\begin{figure}
    \centering
    \includegraphics[width=0.55\linewidth]{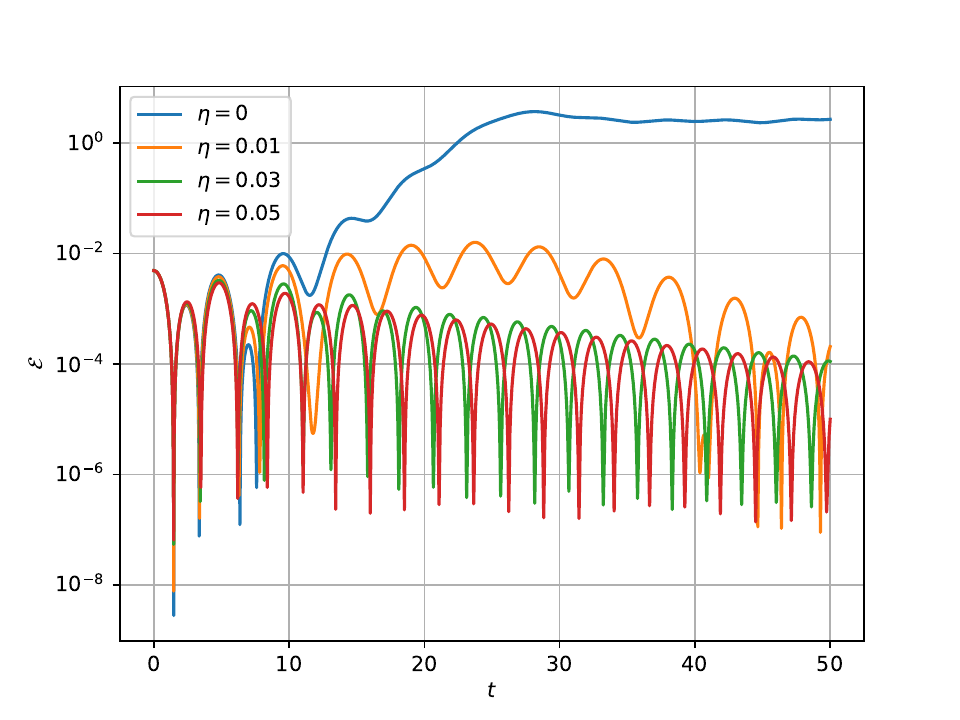}
    \caption{Two-stream instability. Evolution of electric energy for different collision strengths.}
    \label{fig_Vlasov_TSI}
\end{figure}

Another benchmark test for the VAFP equation is the two stream instability. The initial condition is given by
\begin{equation}
    f_0(x,\bv) = \frac{1}{2(2\pi)^{3/2}} (1+A\cos(\kappa x))\left(
        \e^{-(v^{(1)}-v^*)^2/2}
        + \e^{-(v^{(1)}+v^*)^2/2}
    \right)
    \e^{-(v^{(2)})^2/2}
    \e^{-(v^{(3)})^2/2}.
\end{equation}
Follow the parameters in \cite{wang2026dynamical}, we choose $A = 0.005$, $\kappa=0.2$ and $v^*=2.4$.
The spatial domain is $x\in[0,2\pi/\kappa]=[0,10\pi]$
with $N_x = 128$ and periodic boundary condition. The velocity domain is truncated to $[v_{\min},v_{\max}]^3=[-9,9]^3$ with $N_v=128$. The time step is chosen as in (\ref{eq_vlasov_time_step}). 
The prechosen rank is $(r_1,r_2) =(8,8)$.
The electric energy is plotted in Figure \ref{fig_Vlasov_TSI}.
In the collisionless case, the evolution of the electric energy is consistent with the results reported in \cite{wang2026dynamical}. 
For $\eta = 0.03$ and $\eta = 0.05$, the electric energy exhibits a damping behavior similar to that observed in the Landau damping example. 
This is consistent with physical intuition that as the collision strength increases, the solution is driven faster to the Maxwellian equilibrium.

\begin{figure}
     \centering
     \begin{subfigure}[b]{0.24\textwidth}
         \centering
         \includegraphics[width=\textwidth]{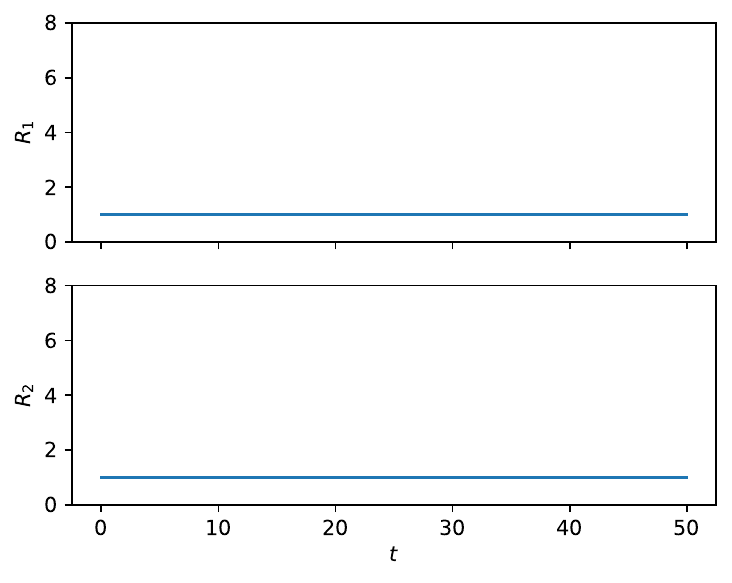}
         \caption{$\eta =0$}
     \end{subfigure}
     \begin{subfigure}[b]{0.24\textwidth}
         \centering
         \includegraphics[width=\textwidth]{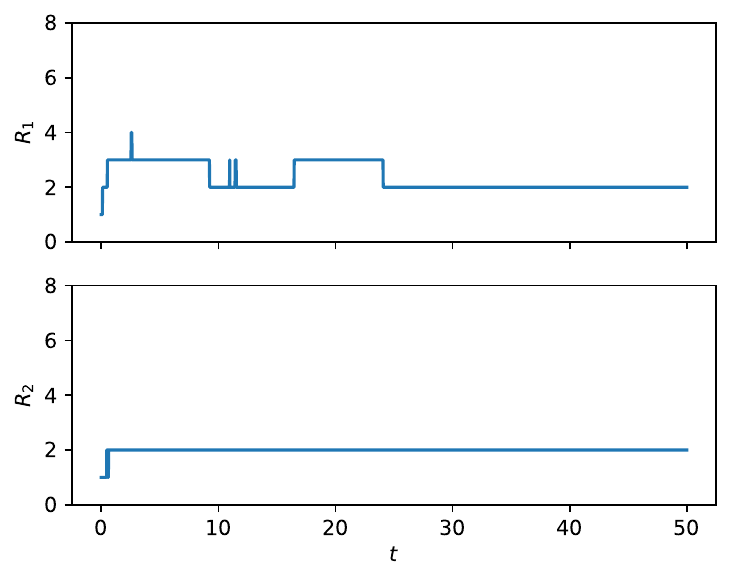}
         \caption{$\eta = 0.01$}
     \end{subfigure}
     \begin{subfigure}[b]{0.24\textwidth}
         \centering
         \includegraphics[width=\textwidth]{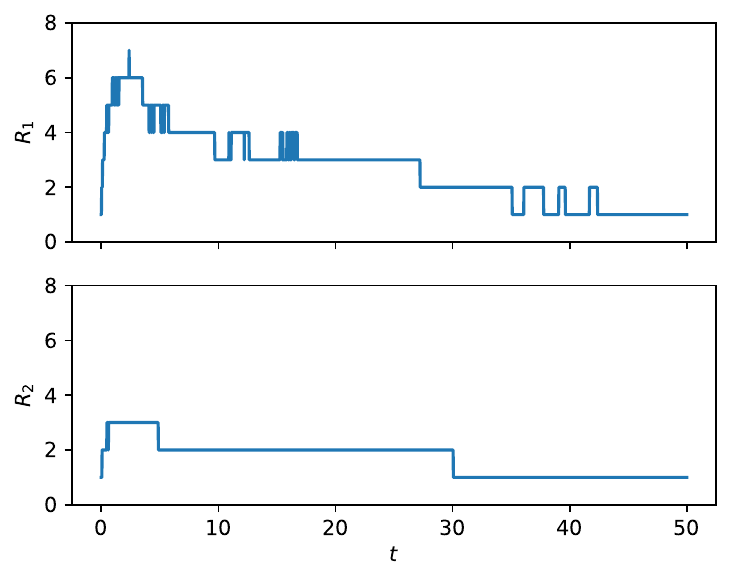}
         \caption{$\eta =0.03$}
     \end{subfigure}
     \begin{subfigure}[b]{0.24\textwidth}
         \centering
         \includegraphics[width=\textwidth]{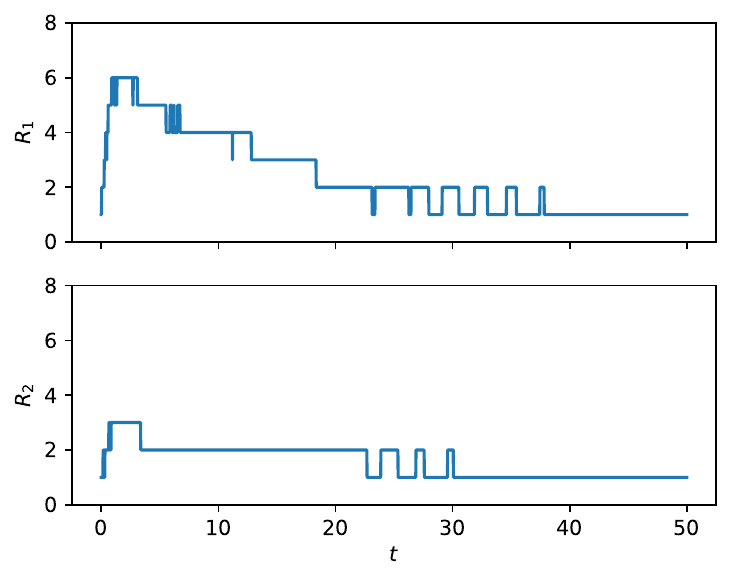}
         \caption{$\eta = 0.05$}
     \end{subfigure}
        \caption{Two-stream instability. Evolution of effective ranks for different collision strengths.}
        \label{fig_TSI_rank}
\end{figure}

\begin{figure}
     \centering
     \begin{subfigure}[b]{0.22\textwidth}
         \centering
         \includegraphics[width=\textwidth]{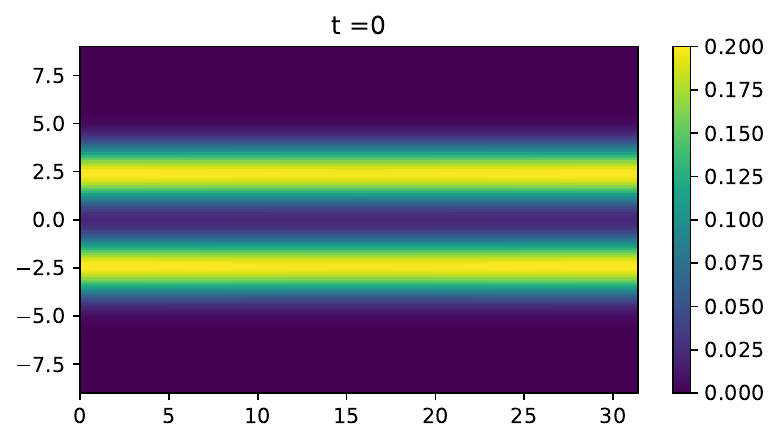}
         \caption{$\eta=0,t=0$}
     \end{subfigure}
     \begin{subfigure}[b]{0.22\textwidth}
         \centering
         \includegraphics[width=\textwidth]{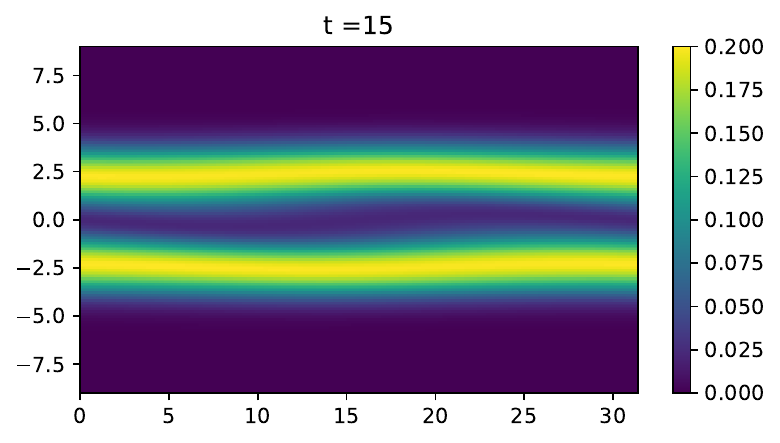}
         \caption{$\eta =0,t=15$}
     \end{subfigure}
     \begin{subfigure}[b]{0.22\textwidth}
         \centering
         \includegraphics[width=\textwidth]{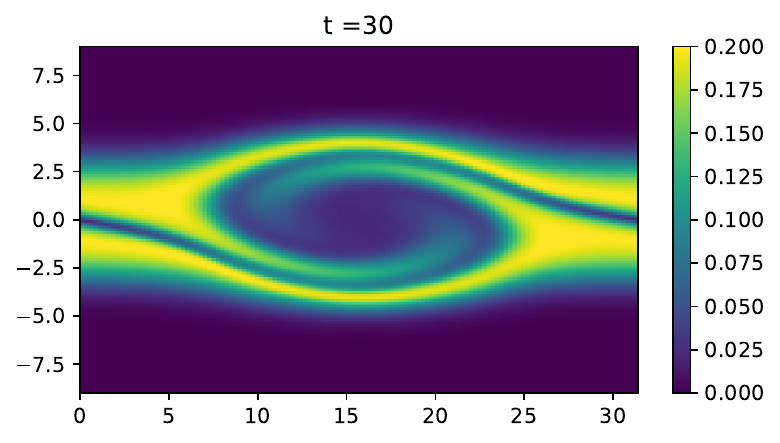}
         \caption{$\eta = 0,t=30$}
     \end{subfigure}
     \begin{subfigure}[b]{0.22\textwidth}
         \centering
         \includegraphics[width=\textwidth]{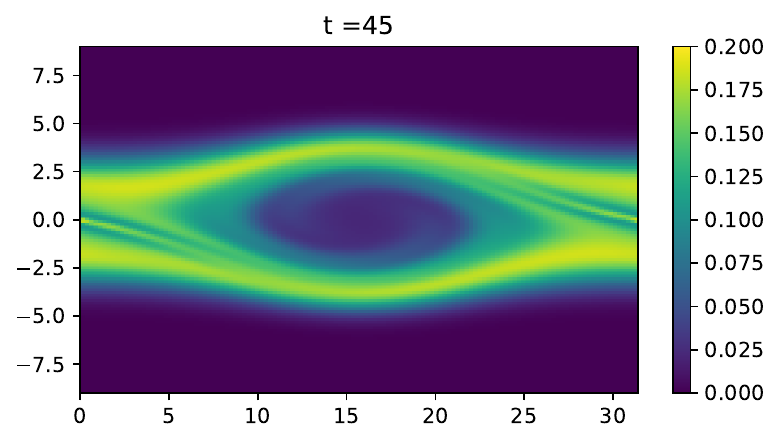}
         \caption{$\eta =0,t=45$}
     \end{subfigure} \\
     \begin{subfigure}[b]{0.22\textwidth}
         \centering
         \includegraphics[width=\textwidth]{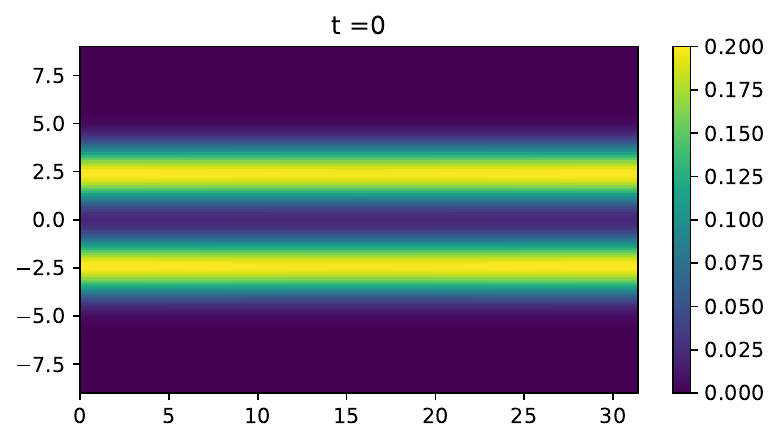}
         \caption{$\eta= 0.01,t=0$}
     \end{subfigure}
     \begin{subfigure}[b]{0.22\textwidth}
         \centering
         \includegraphics[width=\textwidth]{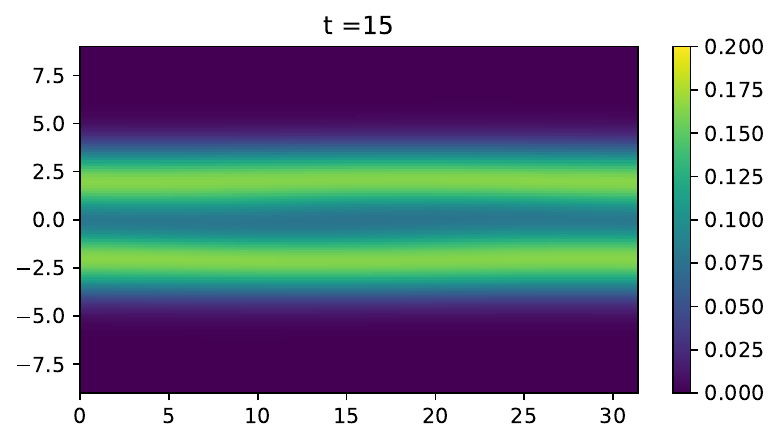}
         \caption{$\eta = 0.01,t=15$}
     \end{subfigure}
     \begin{subfigure}[b]{0.22\textwidth}
         \centering
         \includegraphics[width=\textwidth]{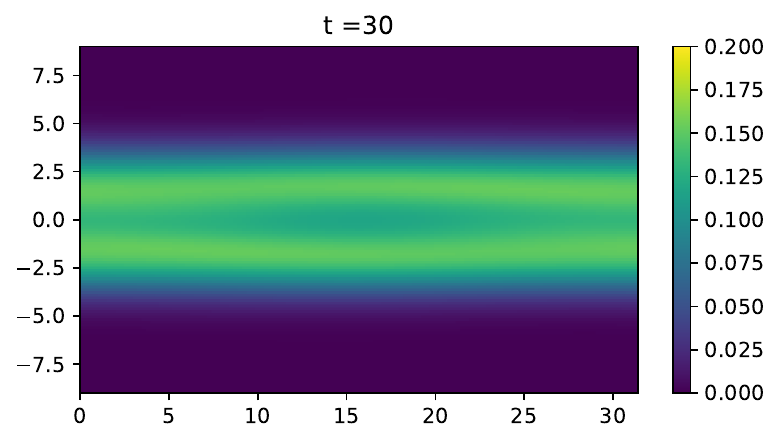}
         \caption{$\eta = 0.01,t=30$}
     \end{subfigure}
     \begin{subfigure}[b]{0.22\textwidth}
         \centering
         \includegraphics[width=\textwidth]{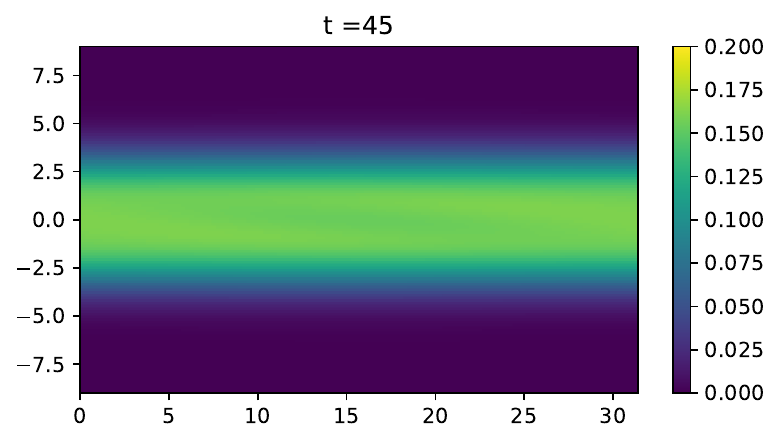}
         \caption{$\eta =0.01,t=45$}
     \end{subfigure}\\
     \begin{subfigure}[b]{0.22\textwidth}
         \centering
         \includegraphics[width=\textwidth]{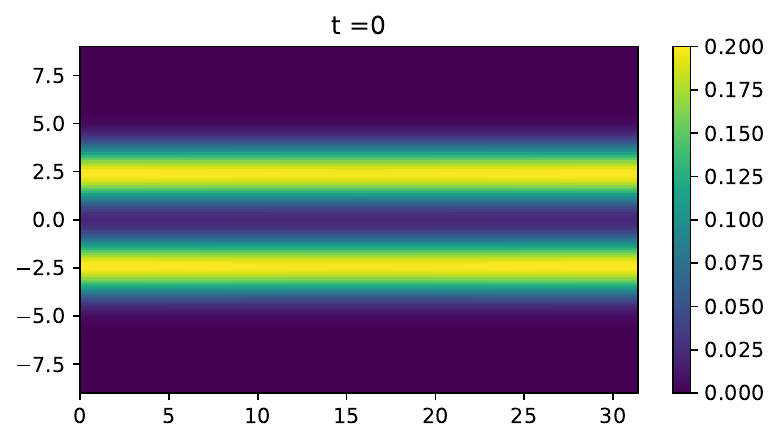}
         \caption{$\eta= 0.03,t=0$}
     \end{subfigure}
     \begin{subfigure}[b]{0.22\textwidth}
         \centering
         \includegraphics[width=\textwidth]{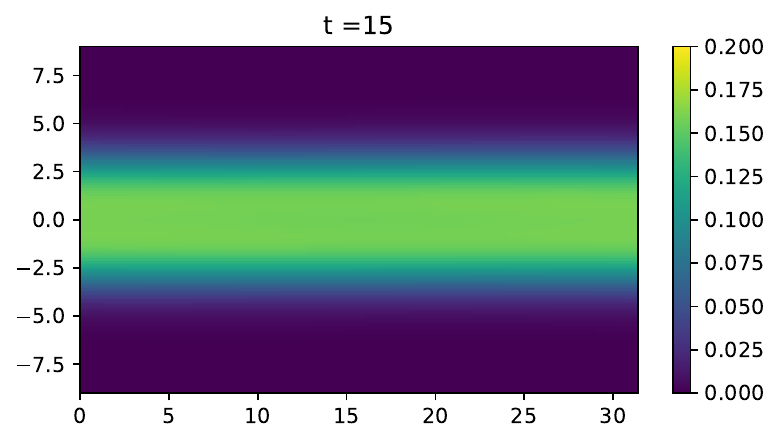}
         \caption{$\eta = 0.03,t=15$}
     \end{subfigure}
     \begin{subfigure}[b]{0.22\textwidth}
         \centering
         \includegraphics[width=\textwidth]{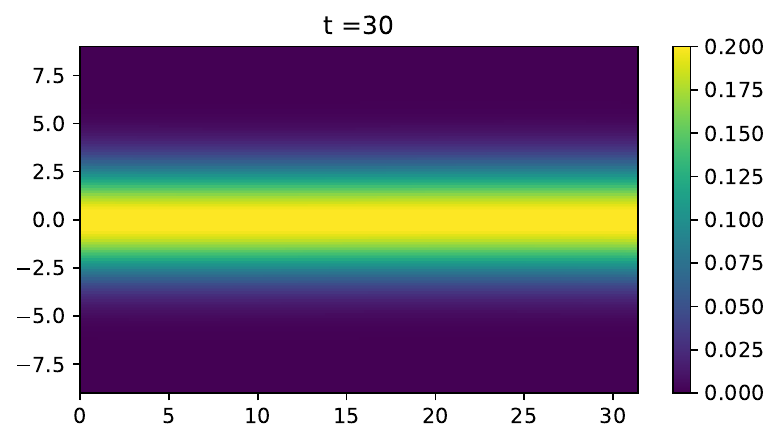}
         \caption{$\eta = 0.03,t=30$}
     \end{subfigure}
     \begin{subfigure}[b]{0.22\textwidth}
         \centering
         \includegraphics[width=\textwidth]{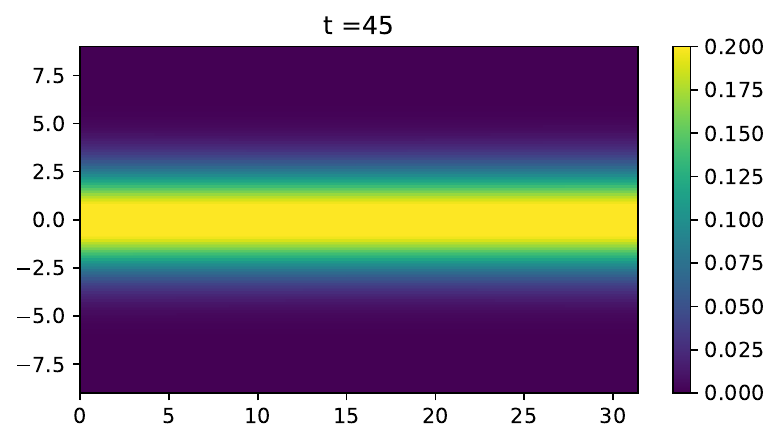}
         \caption{$\eta =0.03,t=45$}
     \end{subfigure}\\
     \begin{subfigure}[b]{0.22\textwidth}
         \centering
         \includegraphics[width=\textwidth]{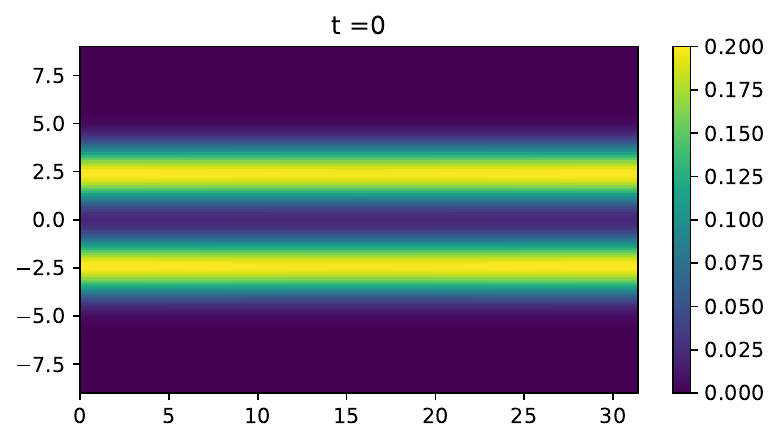}
         \caption{$\eta= 0.05,t=0$}
     \end{subfigure}
     \begin{subfigure}[b]{0.22\textwidth}
         \centering
         \includegraphics[width=\textwidth]{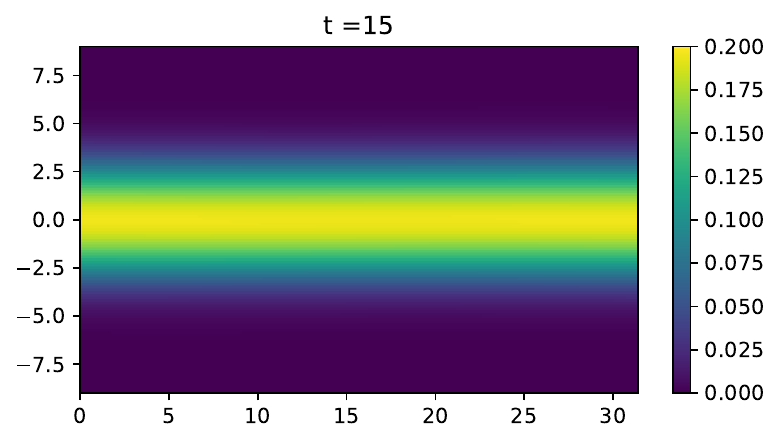}
         \caption{$\eta = 0.05,t=15$}
     \end{subfigure}
     \begin{subfigure}[b]{0.22\textwidth}
         \centering
         \includegraphics[width=\textwidth]{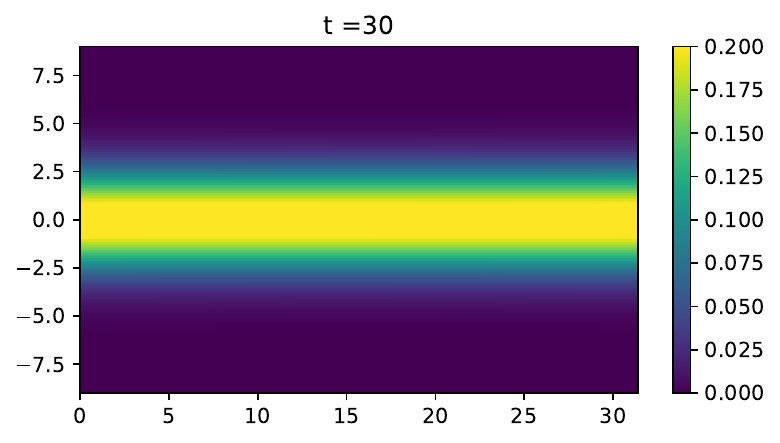}
         \caption{$\eta = 0.05,t=30$}
     \end{subfigure}
     \begin{subfigure}[b]{0.22\textwidth}
         \centering
         \includegraphics[width=\textwidth]{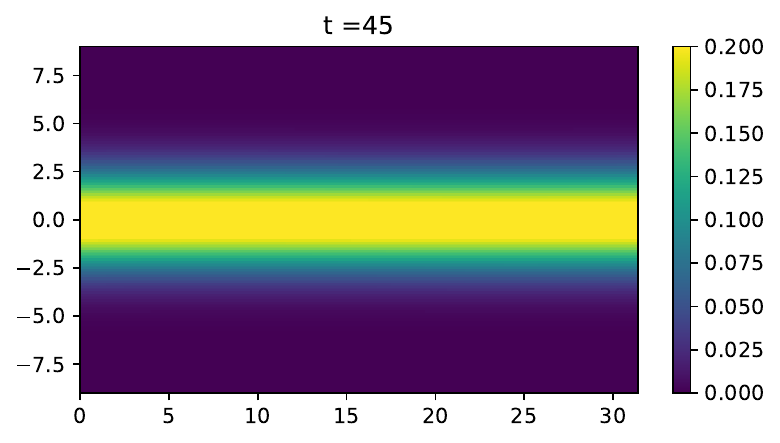}
         \caption{$\eta =0.05,t=45$}
     \end{subfigure}
    \caption{Two-stream instability. Phase plots for different collision strengths.}
    \label{fig_TSI_phase_plot}
\end{figure}

To ensure that the choice of TT-rank is sufficient, we compute the time evolution of the effective rank of the solution.
For a tensor train $\bf$, we compute the singular values $\sigma_1,\cdots,\sigma_{r_1}$ of $S^{(1)}$ when $\bf$ is in form \rombracket{2} of the tensor diagram \eqref{eq_TT_decomposition_forms1}.
The effective rank is then defined as
\begin{equation}
    \mathfrak{r}_1(\bf) = \max \{\ell: \  \sigma_\ell\geq \delta \sigma_1\}.
\end{equation}
Similarly, when $\bf$ is in form \rombracket{4}, another effective rank $\mathfrak{r}_2$ can be defined.
The effective ranks of the solution are then given by $$R_1 = \max_{j=1,\cdots,N_x} \mathfrak{r}_1 (\bf^j), \quad
R_2 = \max_{j=1,\cdots,N_x} \mathfrak{r}_2 (\bf^j).$$
The evolution of the effective ranks with threshold $\delta=10^{-5}$ is shown in Figure \ref{fig_TSI_rank}.
Additionally, the evolution of the phase plots is presented in Figure \ref{fig_TSI_phase_plot}.
In Figure \ref{fig_TSI_rank}, we observe that the effective rank first increases and then decreases in the collisional cases. 
For stronger collisions ($\eta = 0.03$ and $\eta = 0.05$), the effective rank eventually decreases to $(1,1)$. 
This indicates that the solution is numerically driven toward a local Maxwellian state.
The same phenomenon can also be observed in the phase plots in Figure \ref{fig_TSI_phase_plot}. 
In the collisionless case, a clear vortex structure is observed. 
In contrast, in the collisional cases, the vortex gradually smears out and disappears, while the Maxwellian structure emerges.

\section{Conclusion}
\label{sec_conclusion}
In this paper, we present an implicit dynamical tensor-train method for kinetic equations with stiff Fokker--Planck collisions. 
The implicit treatment of the Fokker--Planck operator enables the method to be applied over a wide range of collisional regimes without being restricted by excessively small time steps. In the proposed method, the spatial variables are treated as parameters, while the velocity space is discretized using tensor trains. 
A five-substep procedure is employed, in which three forward substeps are treated implicitly and two backward substeps are handled explicitly. 
In the implicit substeps, Sylvester equations are solved to update the tensor cores. 
By exploiting the special structure of the Fokker--Planck collision operator, we develop tailored matrix and tensor Sylvester solvers.
As a result, the computational cost scales linearly with respect to $N_v$, the number of grid points in a single velocity direction, and remains on the same order as that of an explicit dynamical low-rank method. Several numerical examples are presented to demonstrate the accuracy and efficiency of the proposed method. Future work includes the extension of the method to more general collision operators, such as anisotropic diffusions.

\section*{Acknowledgement}
This work was partially supported by DoD MURI grant FA9550-24-1-0254. The work of JH was additionally supported by NSF grants DMS-2409858 and IIS-2433957, and DOE grant DE-SC0023164.

\bibliographystyle{abbrv}
\bibliography{myBib_abbr}

\end{document}